# SEMI-CLASSICAL ANALYSIS OF A RANDOM WALK ON A MANIFOLD

By Gilles Lebeau and Laurent Michel

*Université de Nice Sophia-Antipolis*

We prove a sharp rate of convergence to stationarity for a natural random walk on a compact Riemannian manifold $(M, g)$. The proof includes a detailed study of the spectral theory of the associated operator.

**1. Introduction.** This paper has two main aims. First, we study the spectral theory of a Markov chain associated to a natural "ball walk" on a compact, connected Riemannian manifold. From $x$, the walk moves to a uniformly chosen point in a ball of radius $h$ around $x$. Here $h$ is a small parameter. We prove a precise Weyl-type estimate on the number of eigenvalues close to 1, and convergence of the spectrum near 1 (when $h \to 0$) to the Laplace–Beltrami spectrum. This walk does not have, in general, the Riemannian area distribution as stationary distribution. The second aim is to analyse the Metropolis algorithm as a way to achieve uniformity. Sharp rates of convergence for the Metropolized chain are given. In the Appendix, we prove that under appropriate scaling, the modified Metropolis chain converges to the Brownian motion.

Let $(M, g)$ be a smooth, compact, connected Riemannian manifold of dimension $d$, equipped with its canonical volume form $d_g x$. Let $d_g(x, y)$ be the Riemannian distance on $M \times M$. For $x \in M$ and $h > 0$, let $B(x, h) = \{y, d_g(x, y) \leq h\}$ be the ball of radius $h$ centered at $x$, and let $|B(x, h)| = \int_{B(x,h)} d_g y$ be its Riemannian volume. For any given $h > 0$, let $T_h$ be the operator acting on continuous functions on $M$,

$$(1.1) \qquad (T_h f)(x) = \frac{1}{|B(x, h)|} \int_{B(x, h)} f(y) \, d_g y.$$

We denote by $K_h$ the kernel of $T_h$, which is given by

$$(1.2) \qquad K_h(x,y)\, d_g y = \frac{\mathbf{1}_{\{d_g(x,y)\leq h\}}}{|B(x,h)|}\, d_g y.$$

Obviously, for any $x \in M$, $K_h(x,y)\, d_g y$ is a probability measure on $M$, and therefore $K_h$ is a Markov kernel. It is associated with the following natural random walk on $M$: if the walk is at $x$, then it moves to a point $y \in B(x,h)$ with a probability given by $K_h(x,y)\, d_g y$.

Let ${}^t T_h$ be the transpose operator acting on Borel measures on $M$, defined as usual by $\langle {}^t T_h(\mu), f\rangle = \langle \mu, T_h(f)\rangle$. Let $c_d$ be the volume of the unit ball of the Euclidean space $\mathbb{R}^d$. For $h$ small, $h^{-d}|B(x,h)|$ is a smooth function on $M$ which converges uniformly on $M$ to $c_d$ when $h \to 0$. Let $d\nu_h$ be the probability measure on $M$,

$$(1.3) \qquad d\nu_h = \frac{|B(x,h)|}{Z_h c_d h^d}\, d_g x,$$

where the normalizing constant $Z_h$ is such that $d\nu_h(M) = 1$. Then for $h$ small, $d\nu_h$ is close to $d_g x / Vol(M)$ and $Z_h$ is close to $Vol(M)$. One verifies easily that $T_h$ is self-adjoint on the space $L^2(M, d\nu_h)$, and that ${}^t T_h(d\nu_h) = d\nu_h$.

The first goal of this paper is to analyze the spectral theory of the self-adjoint operator $T_h$ acting on $L^2(M, d\nu_h)$. Let us recall some basic facts. One has $T_h(1) = 1$, and by the Markov property, the norm of $T_h$ acting on the space $L^\infty$ is equal to 1; by self-adjointness, the norm of $T_h$ acting on the space $L^1(M, d\nu_h)$ is equal to 1 and thus the norm of $T_h$ acting on the space $L^2(M, d\nu_h)$ is also equal to 1. Observe that for any given $h > 0$, the operator $T_h$ is compact. Thus the spectrum of $T_h$, $Spec(T_h)$, is a closed subset of $[-1, 1]$ which is discrete in $[-1, 1] \setminus \{0\}$ with 0 as accumulation point, and each $\mu \in Spec(T_h) \setminus \{0\}$ is an eigenvalue of finite multiplicity.

We denote by $\Delta_g$ the (negative) Laplace–Beltrami operator on $(M, g)$, and by $0 = \lambda_0 < \lambda_1 \leq \lambda_2 \leq \cdots \leq \lambda_n \leq \cdots$ the spectrum of the self-adjoint operator $-\Delta_g$ on $L^2(M, d_g x)$. We will denote by $G_d(\xi)$ the following function of $\xi \in \mathbb{R}^d$:

$$(1.4) \qquad G_d(\xi) = \frac{1}{c_d} \int_{|y|\leq 1} e^{iy\xi}\, dy.$$

Up to the factor $\frac{1}{c_d}$, the function $G_d$ is the Fourier transform of the characteristic function of the unit ball in $\mathbb{R}^d$, and depends only on $|\xi|^2$. We shall also use the function $\Gamma_d(s)$ on $[0, \infty[$ defined by

$$(1.5) \qquad G_d(\xi) = \Gamma_d(|\xi|^2).$$

The function $\Gamma_d$ is real analytic, $|\Gamma_d(s)| \leq 1$, and $\lim_{s \to \infty} \Gamma_d(s) = 0$, since $G_d(\xi)$ is the Fourier transform of a compactly supported, real and even $L^1$ function of total mass 1. One has near $s = 0$,

$$\Gamma_d(s) = 1 - \frac{s}{2(d+2)} + \mathcal{O}(s^2). \tag{1.6}$$

Moreover, there exists $\gamma_0 < 1$ such that $\Gamma_d(s) \in [-\gamma_0, 1]$ for all $s$, and one has $\Gamma_d(s) = 1$ iff $s = 0$. To see this point, just observe that if $|G_d(\xi)| = 1$, then one has $G_d(\xi) = e^{i\theta}$ for some real $\theta$, hence $\int_{|y| \leq 1} (e^{iy\xi - i\theta} - 1)\, dy = 0$ which implies $y\xi - \theta \in 2\pi\mathbb{Z}$ for all $|y| \leq 1$, and therefore $\xi = 0$ and $\theta \in 2\pi\mathbb{Z}$.

THEOREM 1. *Let $h_0 > 0$ be small. There exist $\gamma < 1$ such that for any $h \in ]0, h_0]$, one has $\mathrm{Spec}(T_h) \subset [-\gamma, 1]$, and 1 is a simple eigenvalue of $T_h$. Let*

$$0 < \cdots \leq \mu_{k+1}(h) \leq \mu_k(h) \leq \cdots \leq \mu_1(h) < \mu_0(h) = 1 \tag{1.7}$$

*be the decreasing sequence of positive eigenvalues of $T_h$. For any given $L > 0$, there exists $C$ such that for all $h \in ]0, h_0]$ and all $k \leq L$, one has*

$$\left| \frac{1 - \mu_k(h)}{h^2} - \frac{\lambda_k}{2(d+2)} \right| \leq Ch^2. \tag{1.8}$$

*Let $N(a, h)$ be the number of eigenvalues of $T_h$ in the interval $[a, 1]$. For any given $\delta \in ]0, 1[$, there exist $C_{\delta,i}$ independent of $h \in ]0, h_0]$, such that the following holds true:*

*For any $\tau \in [0, (1-\delta)h^{-2}]$, $N(1 - \tau h^2, h)$ satisfies the Weyl law,*

$$\left| N(1 - \tau h^2, h) - (2\pi h)^{-d} \int_{\Gamma_d(|\xi|_x^2) \in [1-\tau h^2, 1]} dx\, d\xi \right| \tag{1.9}$$
$$\leq C_{\delta,1}(1 + \tau)^{(d-1)/2},$$

*where $dx\, d\xi$ is the canonical volume form on the symplectic manifold $T^*M$, and $|\xi|_x$ is the Riemannian length of the co-vector $\xi$ at $x$. In particular, one has*

$$N(1 - \tau h^2, h) \leq C_{\delta,2}(1 + \tau)^{d/2}. \tag{1.10}$$

*Moreover, for any eigenfunction $e_k^h$ of $T_h$ associated with the eigenvalue $\mu_k(h) \in [\delta, 1]$, the following inequality holds true with $\tau_k(h) = h^{-2}(1 - \mu_k(h))$,*

$$\|e_k^h\|_{L^\infty} \leq C_{\delta,3}(1 + \tau_k(h))^{d/4} \|e_k^h\|_{L^2}. \tag{1.11}$$

Let $|\Delta_h|$ be the positive, bounded, self-adjoint operator on $L^2(M, d\nu_h)$ defined by

$$1 - T_h = \frac{h^2}{2(d+2)} |\Delta_h|. \tag{1.12}$$

By (1.8), the two operators $|\Delta_h|$ and $-\Delta_g$ have almost the same eigenvalues in any interval $[0, L]$ independent of $h$, for $h$ small enough. Our next result gives more precise information on the difference of their resolvents for $h$ small. Observe that as vector spaces, the two Hilbert spaces $L^2(M, d\nu_h)$ and $L^2(M, d_g x)$ are equal, and that their norms are uniformly in $h$ equivalent. We set $L^2 = L^2(M, d\nu_h) = L^2(M, d_g x)$, $\|f\|_{L^2} = \|f\|_{L^2(M, d_g x/Vol(M))}$, and if $A$ is a bounded operator on $L^2$, we denote by $\|A\|_{L^2}$ its norm.

Let $F_1$ and $F_2$ be the two closed subsets of $\mathbb{C}$, $F_1 = \{z, \operatorname{dist}(z, spec(-\Delta_g)) \leq \varepsilon\}$, $F_2 = \{z, \operatorname{Re}(z) \geq A, |\operatorname{Im}(z)| \leq \varepsilon \operatorname{Re}(z)\}$ with $\varepsilon > 0$ small and $A > 0$ large. Let $F = F_1 \cup F_2$ and $U = \mathbb{C} \setminus F$.

THEOREM 2. *There exists $C, h_0 > 0$ such that for all $h \in ]0, h_0]$, and all $z \in U$,*

$$(1.13) \qquad \|(z - |\Delta_h|)^{-1} - (z + \Delta_g)^{-1}\|_{L^2} \leq C h^2.$$

REMARK 1. The error term $\mathcal{O}(h^2)$ in the estimate (1.13) is of the same type than the error one gets for the difference between discrete and continuous Laplacian on $\mathbb{R}^d$. However, in our geometric setting, the Ricci curvature of $M$ contributes also to the error term (see Lemma 3 below), and to get a true discrete Laplacian on the manifold $M$, one will have to discretize the integration process in formula (1.1). Although this is clearly a question of practical interest [as well as modification of $|\Delta_h|$ to improve the convergence in (1.13)], we will not discuss this point in the present paper.

Observe that when $M = (\mathbb{R}/2\pi\mathbb{Z})^d$ is the flat $d$-dimensional torus with $g$ equal to the Euclidean metric, one has the equality,

$$(1.14) \qquad T_h = \Gamma_d(-h^2 \Delta_g).$$

Thus, in that case, the operators $T_h$ and $\Delta_g$ have exactly the same eigenvectors $e^{ikx}$, and the results of Theorems 1 and 2 can be proved by a simple computational verification. For a general compact Riemannian manifold $(M, g)$, the two operators $T_h$ and $\Delta_g$ do not commute, and the formula (1.14) is untrue. In Section 2, we will use a suitable $h$-pseudo-differential calculus in order to show that formula (1.14) remains almost true (in a proper sense), modulo lower order terms involving the curvature of $M$. Then, using the results of Section 2, we will prove Theorems 1 and 2 in Section 3. Observe that the $L^\infty$ bound (1.11) on the eigenfunctions of $|\Delta_h|$ is the exact analogue of what one gets from Sobolev inequalities for the eigenfunctions of $\Delta_g$; in particular, this is certainly not optimal, and it will be of interest to know if the Sogge estimates (see [14]) for the eigenfunctions of $\Delta_g$ are true for the eigenfunctions of $|\Delta_h|$. However, (1.11) will be sufficient for us in the proof of Theorem 3.

Let us now discuss the second goal of this paper. For any $n \geq 1$, let $K_h^n(x,y)\,d_g y$ be the kernel of $(T_h)^n$. Then $\int_A K_h^n(x,y)\,d_g y$ is the probability that the random walk associated to $T_h$ starting at $x$ is in the set $A$ after $n$ steps of the walk. When $n \to \infty$, the sequence of probabilities $K_h^n(x,y)\,d_g y$ will converge to the stationary probability $d\nu_h(y)$, but this is not quite satisfactory, since on a general manifold $M$, $d\nu_h(y)$ depends on $h$. Thus, in order to get a Markov chain with the fixed stationary probability $d\mu_M = d_g x/Vol(M)$, we modified the kernel $K_h(x,y)\,d_g y$, according to the strategy of the Metropolis algorithm, in the following way. Let

(1.15) $$M_h(x,dy) = m_h(x)\delta_{y=x} + \mathcal{K}_h(x,y)\,d_g y,$$

where the functions $m_h$ and $\mathcal{K}_h$ are defined by

(1.16)
$$\mathcal{K}_h(x,y) = K_h(x,y) \min\left(\frac{|B(x,h)|}{|B(y,h)|}, 1\right),$$

$$m_h(x) = 1 - \int_M \mathcal{K}_h(x,y)\,d_g y.$$

Then, $M_h(x,dy)$ is still a Markov kernel, but now, the operator

(1.17) $$M_h(f)(x) = \int_M f(y) M_h(x,dy)$$

is self-adjoint on the space $L^2(M, d_g x)$, and therefore one has ${}^t M_h(d_g x) = d_g x$ for all $h$. Let $M_h^n(x,dy)$ be the kernel of $(M_h)^n$. Our purpose is to get an estimate uniform with respect to the small parameter $h$, on the speed of convergence, when $n \to \infty$, of the probability $M_h^n(x,dy)$ toward the invariant measure $d\mu_M = d_g x/Vol(M)$. Let us recall that if $p,q$ are two probabilities, their total variation distance is defined by

$$\|p - q\|_{TV} = \sup_A |p(A) - q(A)|,$$

where the sup is over all Borel sets $A$. The following theorem tells us that this speed of convergence is estimated for $h$ small, as expected, by the first nonzero eigenvalue $\lambda_1$ of the Laplace–Beltrami operator-$\Delta_g$.

THEOREM 3. *Let $h_0 > 0$ small. There exists $A$ such that for all $h \in ]0,h_0]$ the following holds true:*

$$e^{-\gamma'(h)nh^2} \leq 2 \sup_{x \in M} \|M_h^n(x,dy) - d\mu_M\|_{TV},$$

(1.18)
$$\sup_{x \in M} \|M_h^n(x,dy) - d\mu_M\|_{TV} \leq A e^{-\gamma h nh^2} \quad \text{for all } n.$$

*Here $\gamma(h), \gamma'(h)$ are two positive functions such that $\gamma(h) \simeq \gamma'(h) \simeq \frac{\lambda_1}{2(d+2)}$ when $h \to 0$.*

Of course, the analogue of this result is also valid if one replaces $M_h$ by $T_h$ and $d\mu_M$ by $d\nu_h$, with a simple proof. Theorem 3 will be proved in Section 4. We will verify that $M_h$ is a sufficiently small perturbation of $T_h$, and, in particular, that estimates (1.11) and (1.10) remains true for its eigenfunctions. Finally, in Theorem 4 of the Appendix, we will answer a question of one of the referees of the paper, about the convergence of the Metropolis chain to the Brownian motion on the Riemannian manifold $(M, g)$.

Perhaps the main contribution of this paper is the introduction of microlocal analysis as a tool for analyzing rates of convergence for Markov chains. These result in a fairly general picture; the top of the spectrum of the Metropolis chain converges to a Laplace spectrum. Because of the holding, the Metropolis chain has a continuous spectrum but this is bound from $\pm 1$ and does not enter the final result. This picture was found in a simple case in [4] and for the Metropolis algorithm in Lipschitz domains, including the random placement of $N$ hard discs in the unit square, in [5]. The present paper shows that the picture holds fairly generally. Throughout this paper, we will use basic techniques in semi-classical analysis, for which we refer to [13] and [7].

For an introduction to the well-developed area of probability theory on Riemannian manifolds we refer to [11]. For the analysis of the Metropolis algorithm, we refer to [6] and references therein. There are also emerging applications to statistics on Riemannian manifolds (see [1–3, 10] for examples and references). All of these applications lead to the problem of drawing random samples from the uniform distribution. This topic has not been widely addressed. Some algorithms are suggested in [3]. The present paper is a contribution to a rigorous treatment, giving reasonably sharp bounds on rates of convergence.

**2. The symbolic calculus of $T_h$.** We first recall some basic facts on the classical $h$-pseudo-differential calculus. For $m \in \mathbb{R}$, let $S^m$ the set of functions $a(x, \xi, h)$ smooth in $(x, \xi) \in \mathbb{R}^{2d}$, with parameter $h \in \,]0, 1]$ such that for any $\alpha, \beta$, there exists $C_{\alpha,\beta}$ such that for all $(x, \xi) \in \mathbb{R}^{2d}$ and all $h \in \,]0, 1]$ one has

$$(2.1) \qquad |\partial_x^\alpha \partial_\xi^\beta a(x, \xi, h)| \leq C_{\alpha,\beta}(1 + |\xi|)^{m-|\beta|}.$$

For $a \in S^m$, we denote by $Op(a)$ the $h$-pseudo-differential operator acting on the Schwartz space $\mathcal{S}(\mathbb{R}^d)$,

$$(2.2) \qquad Op(a)(f)(x) = (2\pi h)^{-d} \int e^{i(x-y)\xi/h} a(x, \xi, h) f(y)\, dy\, d\xi.$$

Let us recall that for $a \in S^0$, the operator $Op(a)$ is uniformly bounded in $h$ on the space $L^2(\mathbb{R}^d)$, and that for $a \in S^m, b \in S^k$, one has $Op(a)Op(b) = Op(c)$ where $c = a \sharp b \in S^{m+k}$ is given by the oscillatory integral

(2.3) $\quad c(x, \xi, h) = (2\pi h)^{-d} \int e^{-iz\theta/h} a(x, \xi + \theta, h) b(x + z, \xi, h) \, dz \, d\theta,$

and admits the asymptotic expansion

(2.4)
$$c(x, \xi, h) = \sum_{|\alpha| < N} \frac{h^{|\alpha|}}{i^{|\alpha|} \alpha!} \partial_\xi^\alpha a(x, \xi, h) \partial_x^\alpha b(x, \xi, h)$$
$$+ h^N r_N(x, \xi, h), \qquad r_N \in S^{m+l-N}.$$

The subset $S_{cl}^m$ of $S^m$ is the set of $a(x, \xi, h) \in S^m$ such that there exists a sequence $a_n(x, \xi) \in S^{m-n}, n \geq 0$, such that for all $N$, one has

(2.5) $\quad a(x, \xi, h) = \sum_{0 \leq n < N} (h/i)^n a_n(x, \xi) + h^N r_N(x, \xi, h), \qquad r_n \in S^{m-N}.$

From (2.4), one has $a \sharp b \in S_{cl}^{m+k}$ for $a \in S_{cl}^m$ and $b \in S_{cl}^k$.

Let $(M, g)$ be a compact smooth Riemannian manifold, and let $e_j(x) \in C^\infty(M), j \geq 0$, be an orthonormal basis in $L^2(M, d_g x)$ of real eigenvectors of $-\Delta_g$ with $-\Delta_g e_j = \lambda_j e_j$. For any distribution $f \in \mathcal{D}'(M)$, the Fourier coefficients of $f$ are defined by $f_j = \int f e_j \, d_g x$ and one has $f(x) = \sum_j f_j e_j(x)$ where the series is convergent in $\mathcal{D}'(M)$. For $s \in \mathbb{R}$, let $H^s(M) = (1 - \Delta_g)^{-s/2} L^2(M, d_g x)$ be the usual Sobolev space on $M$. For $f \in \mathcal{D}'(M)$ one has $f \in H^s(M)$ iff $\|f\|_{H^s(M)}^2 = \sum_j (1 + \lambda_j)^s |f_j|^2 < \infty$. We shall also use the semi-classical $H^s$ norms defined by

(2.6) $$\|f\|_{h,s}^2 = \sum_j (1 + h^2 \lambda_j)^s |f_j|^2.$$

A family of operators $R_h$, $h \in ]0, 1]$, acting on the space of distributions $\mathcal{D}'(M)$ is said to be smoothing iff for any $s, t, N$, $R_h$ maps $H^s(M)$ in $H^t(M)$ and there exists $C_{s,t,N}$ such that for all $h \in ]0, 1]$ one has

(2.7) $$\|R_h(f)\|_{H^t(M)} \leq C_{s,t,N} h^N \|R_h(f)\|_{H^s(M)}.$$

A family of operators $A_h$, $h \in ]0, 1]$ acting on the space of distributions $\mathcal{D}'(M)$, belongs to the set $\mathcal{E}_{cl}^m$ of classical $h$-pseudo-differential operators of order $m$, iff for any $x_0 \in M$, there exists an open chart $U$ centered at $x_0$ and two functions $\varphi, \psi \in C_0^\infty(U)$ equal to 1 near $x_0$ with $\psi$ equal to 1 near the support of $\varphi$ such that $A_h \varphi = \psi A_h \varphi + R_h$, with $R_h$ smoothing and there exists $a \simeq \sum_{n \geq 0} (h/i)^n a_n(x, \xi) \in S_{cl}^m$, such that in the local chart $U$, one has $\psi A_h \varphi = Op(a)$. The principal symbol of $A_h$, $\sigma_0(A_h)(x, \xi)$, is by definition the first term $a_0(x, \xi)$ in the asymptotic expansion of $a(x, \xi, h)$. It is a well-defined function on $T^*M$, and for any smooth function $\varphi \in C^\infty(M)$, one has

$$(2.8) \quad e^{-i\varphi(x)/h} A_h(e^{i\varphi(x)/h}) = \sigma_0(A_h)(x, d\varphi(x)) + \mathcal{O}(h).$$

Then $\mathcal{E}_{cl} = \bigcup_m \mathcal{E}_{cl}^m$ is the algebra of classical $h$-pseudo-differential operators on $M$. For $A_h \in \mathcal{E}_{cl}^m$ and $B_h \in \mathcal{E}_{cl}^k$, one has $A_h B_h \in \mathcal{E}_{cl}^{m+k}$, $\sigma_0(A_h B_h) = \sigma_0(A_h)\sigma_0(B_h)$ and the commutator $[A_h, B_h] = A_h B_h - B_h A_h$ satisfies $[A_h, B_h] \in h\mathcal{E}_{cl}^{m+k-1}$, $\sigma_0(\frac{i}{h}[A_h, B_h]) = \{\sigma_0(A_h), \sigma_0(B_h)\}$ where $\{f, g\}$ is the Poisson bracket. Moreover, for any $A_h \in \mathcal{E}_{cl}^m$, one has $A_h^* \in \mathcal{E}_{cl}^m$, $\sigma_0(A_h^*) = \overline{\sigma_0(A_h)}$, and for any $s \in \mathbb{R}$, there exist $C_s$ independent of $h \in ]0,1]$ such that

$$(2.9) \quad \|A_h f\|_{h,s-m} \leq C_s \|f\|_{h,s} \quad \forall f \in H^s(M).$$

Let us recall that for any $\Phi \in C_0^\infty([0, \infty[)$, the operator $\Phi(-h^2\Delta_g)$ defined by

$$(2.10) \quad \Phi(-h^2\Delta_g)(f) = \sum_j \Phi(h^2\lambda_j) f_j e_j(x)$$

belongs to $\mathcal{E}_{cl}^{-\infty} = \bigcap_m \mathcal{E}_{cl}^m$, and its principal symbol is equal to

$$(2.11) \quad \sigma_0(\Phi(-h^2\Delta_g)) = \Phi(|\xi|_x^2),$$

where $|\xi|_x$ is the Riemannian length of the co-vector $\xi$ at $x$. For a proof of this fact, we refer to [7].

DEFINITION 1. A family of operators $C_h$, $h \in ]0,1]$, acting on the space of distributions $\mathcal{D}'(M)$, belongs to the class $\widetilde{\mathcal{E}}_{cl}^0$ if and only if $C_h$ is bounded uniformly in $h$ on $L^2(M)$ and for any $\Phi_0 \in C_0^\infty([0, \infty[)$, one has

$$(2.12) \quad \Phi_0(-h^2\Delta_g)C_h \text{ and } C_h \Phi_0(-h^2\Delta_g) \text{ belongs to } \mathcal{E}_{cl}^{-\infty}.$$

Let $\Gamma_{d,h}$ be the operator $\Gamma_{d,h} = \Gamma_d(-h^2\Delta_g)$, so that

$$(2.13) \quad \Gamma_{d,h}(f)(x) = \sum_j \Gamma_d(h^2\lambda_j) f_j e_j(x).$$

Since $\Phi_0\Gamma_d \in C_0^\infty([0, \infty[)$, one has obviously $\Gamma_{d,h} \in \widetilde{\mathcal{E}}_{cl}^0$.

Let $U \subset M$ be an open chart with local coordinates $x = (x_1, \ldots, x_d) \in \mathbb{R}^d$. Then for $x \in U$ and $r > 0$ small, the geodesic ball of radius $r$ centered at $x$ is given by

$$(2.14) \quad B(x, r) = \Big\{x + u, \sum k_{i,j}(x, u) u_i u_j \leq r^2\Big\},$$

where $(k_{i,j}(x,u))$ is a smooth and symmetric matrix in $(x,u)$ such that $k_{i,j}(x,0) = g_{i,j}(x)$. For any function $f$ compactly supported in $U$ and $h$ small, $T_h f$ is supported in $U$ and given in these local coordinates by

$$(2.15) \quad T_h f(x) = \frac{1}{|B(x,h)|} \int_{{}^t u k(x,u) u \leq h^2} f(x+u) \sqrt{\det(g(x+u))} \, du.$$

Using the new integration variable $hv = w = k^{1/2}(x,u)u$ in (2.15), we get

$$(2.16) \quad T_h f(x) = \frac{h^d}{|B(x,h)|} \int_{|v| \leq 1} f(x + hm(x,hv)v) \rho(x,hv) \, dv,$$

where $m(x,w)$ is the smooth, symmetric and positive matrix, such that near $u = 0$ one has $w = k^{1/2}(x,u)u \Leftrightarrow u = m(x,w)w$, so $m(x,0) = g^{-1/2}(x)$, and $\rho(x,w) = \sqrt{\det(g(x+u))} |\det \frac{\partial u}{\partial w}|$ is smooth in $(x,w)$ and $\rho(x,0) = 1$.

LEMMA 1. *For $h_0 > 0$ small and any $k$, $T_h$ is a bounded operator on $C^k(M)$ uniformly in $h \in ]0, h_0]$. Moreover, there exists $C$ independent of $h$ such that, with $|\Delta_h|$ defined in (1.12), one has for all $f \in C^2(M)$,*

$$(2.17) \quad \||\Delta_h| f\|_{L^\infty} \leq C \|f\|_{C^2}.$$

PROOF. The first assertion is obvious from (2.16) since $\frac{h^d}{|B(x,h)|}$ is a smooth function of $x, h \in [0, h_0]$. From (2.16) and the Taylor formula $f(x+y) = f(x) + \nabla f(x) y + \mathcal{O}(y^2 \|f\|_{C^2})$, one gets easily that (2.17) holds true. □

In the above open chart $U$, we define the symbol of $T_h$, $\sigma(T_h)$ by

$$(2.18) \quad \sigma(T_h)(x, \xi, h) = e^{-ix\xi/h} T_h(e^{ix\xi/h}).$$

For any compact set $K \subset U$, there exists $h_K > 0$ such that $\sigma(T_h)(x, \xi, h)$ is well defined for $x \in K, \xi \in \mathbb{R}^d$ and $h \in ]0, h_K]$. From (2.18), one has

$$(2.19) \quad \sigma(T_h)(x, \xi, h) = \frac{h^d}{|B(x,h)|} \int_{|v| \leq 1} e^{i{}^t \xi . m(x,hv) v} \rho(x, hv) \, dv,$$

and therefore, for any $\alpha, \beta$, there exists $C_{\alpha,\beta}$ independent of $h$ such that

$$(2.20) \quad |\partial_x^\alpha \partial_\xi^\beta \sigma(T_h)(x, \xi, h)| \leq C_{\alpha,\beta}(1 + |\xi|)^{|\alpha|}.$$

Observe also that, since $m(x,0) = g^{-1/2}(x)$ and $\rho(x,0) = 1$, one has

$$(2.21) \quad \sigma(T_h)(x, \xi, 0) = \Gamma_d(|\xi|_x^2).$$

LEMMA 2. *Let $h_0$ small. For $h \in ]0, h_0]$, the operator $T_h$ belongs to the class $\widetilde{\mathcal{E}}_{cl}^0$.*

PROOF. Let $M = \bigcup_k U_k$ be a finite covering of $M$ by local charts $U_k$, and $1 = \sum_k \varphi_k(x)$ a partition of unity with $\varphi_k \in C_0^\infty(U_k)$. Let $\psi_k \in C_0^\infty(U_k)$ equal to 1 near the support of $\varphi_k$. Then for $h$ small enough, one has

$$(2.22) \qquad T_h(f)(x) = \sum_k \psi_k T_h(\varphi_k f)(x).$$

Let $T_{h,k} = \psi_k T_h \varphi_k$; we reduce to show that for any $k$, $T_{h,k} \in \widetilde{\mathcal{E}}_{cl}^0$. Let $\Phi_0 \in C_0^\infty[0, \infty[$; there exists $\psi \in C_0^\infty(U_k)$ and a compact set $K \subset U_k$ such that $\varphi_k \Phi_0(-h^2 \Delta_g) = Op(a)\psi + R_h$ with $a(x, \xi, h) \in S_{cl}^{-\infty}$ with support in $x \in K$, and $R_h$ smoothing. By Lemma 1, $T_h R_h$ is smoothing, and thus we are reduce to show that in the local chart $U_k$, one has $T_h Op(a) \in S_{cl}^{-\infty}$. From (2.2) and (2.16), one has

$$T_h Op(a)(f)(x) = (2\pi h)^{-d} \int e^{i(x-y)\xi/h} b(x, \xi, h) f(y) \, dy \, d\xi,$$

$$(2.23) \qquad b(x, \xi, h) = \frac{h^d}{|B(x, h|} \int_{|v| \leq 1} e^{it\xi . m(x, hv)v} a(x + hm(x, hv)v, \xi, h)$$

$$\times \rho(x, hv) \, dv.$$

From (2.23) and $a \in S^{-\infty}$, it is clear that $b \in S^{-\infty}$. Using the Taylor expansion in $h$ in (2.23) and $a \in S_{cl}^{-\infty}$, one gets easily $b \in S_{cl}^{-\infty}$. Thus $T_h \Phi_0(-h^2 \Delta_g) \in \mathcal{E}_{cl}^{-\infty}$, and since $T_h$ is self-adjoint for the volume form $d\nu_h$ given by (1.3), one has also $\Phi_0(-h^2 \Delta_g) T_h \in \mathcal{E}_{cl}^{-\infty}$. The proof of our lemma is complete □

Using the Taylor expansion $a(x + hmv, \xi, h) = \sum \frac{(hmv)^\alpha}{\alpha!} \partial_x^\alpha a(x, \xi, h)$ and $(mv)^\alpha e^{it\xi.mv} = (\partial_\xi/i)^\alpha e^{it\xi.mv}$, we get from (2.23) that the symbol $b$ admits the usual asymptotic development,

$$(2.24) \qquad b(x, \xi, h) \simeq \sum_\alpha (h/i)^\alpha \frac{1}{\alpha!} \partial_\xi^\alpha \sigma(T_h)(x, \xi, h) \, \partial_x^\alpha a(x, \xi, h).$$

The following lemma will be crucial in our analysis.

LEMMA 3. *Let $\Phi_0 \in C_0^\infty([0, \infty[)$, and $A_h = h^{-2}(T_h - \Gamma_{d,h})\Phi_0(-h^2 \Delta_g)$. Then $A_h$ belongs to $\mathcal{E}_{cl}^{-\infty}$. Its principal symbol, $\sigma_0(A_h)$, satisfies near $\xi = 0$,*

$$\sigma_0(A_h)(x, \xi) = \left( \frac{S(x)}{3} |\xi|_x^2 (\Gamma_d''(0) - \Gamma_d'(0)^2) \right.$$
$$(2.25)$$
$$\left. + \frac{\Gamma_d''(0)}{3} Ric(x)(\xi, \xi) \right) \Phi_0(|\xi|_x^2) + \mathcal{O}(\xi^3),$$

*where $Ric(x)$ and $S(x)$ are the Ricci tensor and the scalar curvature at $x$. Moreover, let $U$ be a local chart, $K$ a compact subset of $U$ and $\varphi \in C_0^\infty(U)$*

such that $\varphi(x) = 1$ in a neighborhood of $K$; let $a(x, \xi, h) \simeq \sum (h/i)^k a_k(x, \xi) \in S_{cl}^{-\infty}$ be such that in this local chart one has $A_h \varphi = Op(a) + R_h$ with $R_h$ smoothing. Then, for all $k$ and all $x \in K$ one has $a_k(x, 0) = 0$.

PROOF. Let $x_0 \in M$ and let $e_1, \ldots, e_d$ be an orthonormal basis of the tangent space $T_{x_0} M$. For $x = (x_1, \ldots, x_d) \in \mathbb{R}^d$, we identify $x$ with $\sum x_j e_j \in T_{x_0} M$. Let $s \mapsto \exp_{x_0}(sx)$ be the geodesic curve starting at $x_0$ with speed $x$. Then, for $r > 0$ small, the map $\phi_{x_0} : x \mapsto \exp_{x_0}(x)$ is a diffeomorphism of the Euclidean ball $|x| < r$ on an open neighborhood $U$ of $x_0$, and the coordinates $x_j$ in $U$ are called geodesics coordinates centered at $x_0$. In these coordinates, one has $x_0 = 0$, and $(g_{i,j}(0)) = Id$. Let $R$ be the Riemann curvature tensor at $x = 0$ and $R_{(j,k)(l,m)} = (R(\frac{\partial}{\partial x_l}, \frac{\partial}{\partial x_m}) \frac{\partial}{\partial x_k} | \frac{\partial}{\partial x_j})$. Then the Ricci tensor and the scalar curvature at $x = 0$ are given by

$$(2.26) \quad Ric\left(\frac{\partial}{\partial x_j}, \frac{\partial}{\partial x_k}\right) = Ric_{j,k} = \sum_i R_{(i,j)(i,k)}, \qquad S = \sum_j Ric_{j,j}.$$

Moreover, one has in these geodesic coordinates (see [15], page 474)

$$(2.27) \quad \partial_j g_{l,m}(0) = 0, \qquad \partial_j \partial_k g_{l,m}(0) = -\tfrac{1}{3} R_{(l,j)(m,k)} - \tfrac{1}{3} R_{(l,k)(m,j)}$$

or, equivalently,

$$(2.28) \qquad g_{i,j}(x) = \delta_{i,j} + \tfrac{1}{3}(R(x, e_i) x | e_j) + \mathcal{O}(x^3).$$

Consequently, one has

$$(2.29) \qquad \sqrt{\det(g)(x)} = 1 - \tfrac{1}{6} Ric(x, x) + \mathcal{O}(x^3).$$

From this formula, parity arguments, and $2 c_d \Gamma'_d(0) = -\int_{|y| \leq 1} y_j^2 \, dy$, we get

$$(2.30) \qquad |B(0, h)| = h^d c_d \left(1 + \frac{\Gamma'_d(0)}{3} S h^2 + \mathcal{O}(h^3)\right).$$

Moreover, in geodesic coordinates, one has $k(0, u) = Id = m(0, w)$ and $\rho(0, v) = \sqrt{\det(g)(v)}$, and thus from (2.19), (2.29), (2.30) and (1.4), we get

$$\sigma(T_h)(0, \xi, h) = \frac{h^d}{|B(0, h)|} \int_{|v| \leq 1} e^{i\xi \cdot v} \sqrt{\det(g)(hv)} \, dv$$

$$= \Gamma_d(|\xi|^2)\left(1 - \frac{\Gamma'_d(0)}{3} S h^2\right) - \frac{h^2}{6 c_d} \int_{|v| \leq 1} e^{i\xi \cdot v} Ric(v, v) \, dv$$

$$(2.31) \qquad + \mathcal{O}(h^3)$$

$$= \Gamma_d(|\xi|^2) + h^2 \left(-\Gamma_d(|\xi|^2) \frac{\Gamma'_d(0)}{3} S + \frac{1}{6} \sum Ric_{j,k} \frac{\partial^2 G_d}{\partial \xi_j \, \partial \xi_k}(\xi)\right)$$

$$+ \mathcal{O}(h^3).$$

Since $G_d(\xi) = \Gamma_d(|\xi|^2)$, one has

$$\frac{\partial^2 G_d}{\partial \xi_j \, \partial \xi_k}(\xi) = 2\delta_{j,k}(\Gamma'_d(0) + |\xi|^2 \Gamma''_d(0)) + 4\xi_j \xi_k \Gamma''_d(0) + \mathcal{O}(|\xi|^4),$$

and from $\Gamma_d(|\xi|^2) = 1 + \Gamma'_d(0)|\xi|^2 + \mathcal{O}(|\xi|^4)$, we get from (2.31),

$$\begin{aligned}
\sigma(T_h)(0, \xi, h) \\
= \Gamma_d(|\xi|^2) \\
+ h^2 \left( \frac{S|\xi|^2}{3}(\Gamma''_d(0) - \Gamma'_d(0)^2) + \frac{2\Gamma''_d(0)}{3} Ric(\xi, \xi) + \mathcal{O}(|\xi|^4) \right) \\
+ \mathcal{O}(h^3).
\end{aligned} \quad (2.32)$$

Let us now compute the symbol of the operator $\Gamma_{d,h} \Phi_0(-h^2 \Delta_g)$. Until the end of the proof we use the Einstein summation convention. First we remark that in local coordinates the symbol of the operator $-h^2 \Delta_g$ is given by $p = p_0 + h p_1$ with $p_0(x, \xi) = g^{jk}(x) \xi_j \xi_k = |\xi|_x^2$ and $p_1(x, \xi) = -i \tilde{g}_k \xi_k$. Here $(g^{jk})$ denotes the inverse matrix of the matrix $(g_{jk})$ and $\tilde{g}_k = \partial_{x_j} g^{jk} + \frac{1}{2g} g^{jk} \partial_{x_j} g$ where $g$ is the determinant of the matrix $(g_{jk})$. Let $F = \Phi_0 \Gamma_d$ and $\tilde{F}$ be an almost analytic extension of $F$. Then

$$F(-h^2 \Delta_g) = \frac{1}{\pi} \int_{\mathbb{C}} \overline{\partial} \tilde{F}(z)(-h^2 \Delta_g - z)^{-1} L(dz), \quad (2.33)$$

where $L(dz) = dx \, dy$ is the Lebesgue measure on $\mathbb{C}$ and $\overline{\partial} = \frac{1}{2}(\partial_x + i \partial_y)$. Let $\varphi \in C_0^\infty$ be equal to 1 near $x = 0$. For any $z \in \mathbb{C} \setminus \mathbb{R}$ there exist symbols $a_0, a_1, a_2$ such that in local geodesic coordinates we have

$$(-h^2 \Delta_g - z) Op(a_0 + h a_1 + h^2 a_2) = \varphi(x) + h^3 R_h \quad (2.34)$$

with $R_h \in \mathcal{E}_{cl}^0$. From the symbolic calculus it suffices to set

$$\begin{aligned}
a_0 &= \frac{\varphi(x)}{p_0 - z}, & a_1 &= \frac{-i}{p_0 - z} \partial_{\xi_j} p_1 \, \partial_{x_j} a_0, \\
a_2 &= \frac{-1}{p_0 - z}(p_0 \sharp_1 a_1 + p_0 \sharp_2 a_0 + p_1 a_1 + p_1 \sharp_1 a_0),
\end{aligned} \quad (2.35)$$

where for two symbols $f, g$ we define $f \sharp_j g(x, \xi) = \sum_{|\alpha|=j} \frac{1}{i^j \alpha!} \partial_\xi^\alpha f(x, \xi) \partial_x^\alpha g(x, \xi)$. It follows that

$$F(-h^2 \Delta_g) \varphi(x) = Op(b_0 + h b_1 + h^2 b_2) + h^3 \tilde{R}_h \quad (2.36)$$

with $b_j(x, \xi) = \frac{1}{\pi} \int_{\mathbb{C}} \overline{\partial} \tilde{F}(z) a_j(z, x, \xi) L(dz)$ and $\tilde{R}_h \in \mathcal{E}_{cl}^0$. In particular we have $b_0 = \varphi(x) F(|\xi|_x^2)$, and as $a_1(z, 0, \xi) = 0$; we get also $b_1(0, \xi) = 0$. Let us compute $a_2(z, 0, \xi)$. First, we observe that $p_1(0, \xi) = 0$. Moreover, as $\partial_{x_k} p_0(0, \xi) =$ 0, for all $k$ we have also $(p_1 \sharp_1 a_0)(z, 0, \xi) = 0$, $p_0 \sharp_1 a_1(z, 0, \xi) = O(\frac{|\xi|^3}{|\operatorname{Im} z|^3})$ and $p_0 \sharp_2 a_0(z, 0, \xi) = \frac{\Delta g_{lm}(0)\xi_l \xi_m}{(|\xi|^2 - z)^2}$. Therefore, from (2.27) we get

$$
\begin{aligned}
b_2(0, \xi) &= \frac{-1}{\pi} \int_{\mathbb{C}} \overline{\partial} F(z) \frac{1}{(|\xi|^2 - z)^3} L(dz) \Delta g_{lm}(0) \xi_l \xi_m + O(|\xi|^3) \\
&= -\frac{1}{2} F''(|\xi|^2) \Delta g_{lm}(0) \xi_l \xi_m + O(|\xi|^3) \\
&= \frac{1}{3} F''(0) Ric_{lm} \xi_l \xi_m + O(|\xi|^3).
\end{aligned}
\tag{2.37}
$$

Therefore, we conclude that in geodesic coordinates, the symbol of $F(-h^2 \Delta_g)$ satisfies

$$
\begin{aligned}
&\sigma(F(-h^2 \Delta_g))(0, \xi, h) \\
&= F(|\xi|^2) + h^2 \left( \frac{F''(0)}{3} Ric(\xi, \xi) + \mathcal{O}(\xi^3) \right) + \mathcal{O}(h^3).
\end{aligned}
\tag{2.38}
$$

Then, from (2.32), (2.38) and the rule of symbolic calculus, which are valid for $T_h$ by (2.24), we conclude that $A_h$ belongs to $\mathcal{E}_{cl}^{-\infty}$ and that (2.25) holds true.

Finally, since $T_h(1) = 1 = \Gamma_d(-h^2 \Delta_g)(1)$ and $\Phi_0(-h^2 \Delta_g)(1) = \Phi_0(0)$, one has $A_h(1) = 0$; therefore $A_h \varphi(x) = \mathcal{O}(h^\infty)$ for any $x \in K$, and therefore, $Op(a)(1)(x) = a(x, 0, h) = \mathcal{O}(h^\infty)$ for any $x \in K$. The proof of Lemma 3 is complete. □

The following lemma will be used in the sequel to handle the very high frequencies.

LEMMA 4. *Let $\chi \in C_0^\infty(\mathbb{R})$ be equal to 1 near 0. There exists $h_0 > 0, C_0$ such that for all $p \in [1, \infty]$, all $h \in ]0, h_0]$ and all $s \geq 1$, one has*

$$
\left\| T_h(1 - \chi)\left(\frac{-h^2 \Delta_g}{s}\right) \right\|_{L^p} \leq \frac{C_0}{\sqrt{s}}.
\tag{2.39}
$$

PROOF. Set $\hbar = h/\sqrt{s}$. Then $\chi(\frac{-h^2 \Delta_g}{s})$ is a $\hbar$ classical pseudo-differential operator, and belongs to the class $\mathcal{E}_{cl}^{-\infty}$. Let $R_\hbar(x, y) \, d_g y$ be the kernel of the operator $\chi(\frac{-h^2 \Delta_g}{s})$. Then $R_\hbar(x, y)$ is a smooth function of $(x, y) \in M \times M$, and for any $\alpha$, there exists a nonincreasing function $\psi_\alpha$ with rapid decay such that for all $\hbar \in ]0, 1]$, one has

$$
|\nabla_{x,y}^\alpha R_\hbar(x, y)| \leq \hbar^{-d - |\alpha|} \psi_\alpha\left( \frac{d_g(x, y)}{\hbar} \right).
\tag{2.40}
$$

Let $\Theta_{h,s}(x,y)\,d_g y$ be the kernel of the operator $T_h(1-\chi)(\frac{-h^2\Delta_g}{s})$. Then one has

$$\Theta_{h,s}(x,y) = \frac{\mathbf{1}_{\{d_g(x,y)\le h\}}}{|B(x,h)|} - \frac{1}{|B(x,h)|}\int_{B(x,h)} R_\hbar(z,y)\,d_g z. \tag{2.41}$$

By the Shur lemma, it is sufficient to prove that there exists $h_0 > 0, C_0$ such that

$$\sup_{x\in M, h\in ]0,h_0]} \int |\Theta_{h,s}(x,y)|\,d_g y \le C_0/\sqrt{s},$$

$$\sup_{y\in M, h\in ]0,h_0]} \int |\Theta_{h,s}(x,y)|\,d_g x \le C_0/\sqrt{s}. \tag{2.42}$$

We shall prove the first line in (2.42), the proof of the second line being the same. One has $\hbar \le h$ for $s \ge 1$, and from (2.41) and (2.40), we get that for any given $c_0 > 0$, one has for all $h \in ]0, c_0/2]$,

$$\begin{aligned}&d_g(x,y) \ge c_0 \\ &\implies |\Theta_{h,s}(x,y)| \le \hbar^{-d}\psi_0(c_0/2\hbar) \in O(\hbar^\infty) \subset O(s^{-\infty}).\end{aligned} \tag{2.43}$$

Thus we may work in a local chart $U$ centered at a given $x_0 \in M$, with local coordinates $x = (x_1,\ldots,x_d) \in \mathbb{R}^d$, and we are reduced to prove in this local chart, for some $C_0 > 0$ independent of $x_0, h \in ]0,h_0], s \ge 1$,

$$\sup_{h\in ]0,h_0]} \int_{|y|\le 2c_0} |\Theta_{h,s}(x_0=0,y)|\,d_g y \le C_0/\sqrt{s}. \tag{2.44}$$

Let $f_x(y) = \frac{\mathbf{1}_{\{d_g(x,y)\le h\}}}{|B(x,h)|}$. One has

$$\Theta_{h,s}(x_0,y) = f_{x_0}(y) - \int {}^t R_\hbar(y,z) f_{x_0}(z)\,d_g z. \tag{2.45}$$

Let $r_\hbar(y,\xi,\hbar) \simeq \sum_k \hbar^k r_{\hbar,k}(y,\xi) \in S_{cl}^{-\infty}$ be the symbol of ${}^t R_\hbar = \chi(-\hbar^2\Delta_g) \in \mathcal{E}_{cl}^{-\infty}$ in the local chart $U$. Then all the $r_{\hbar,k}(y,\xi)$ are smooth functions of $(y,\xi)$ with support in $|\xi|_y^2 \le r_0$ if $\chi(r)$ is supported in $r \le r_0$. Moreover, by (2.11), one has $r_{\hbar,0}(y,0) = 1$. Therefore, we get with $b_0(y,u)$ smooth in $y$ and in the Schwartz class in $u$, and for some $\psi$ with rapid decay,

$$\begin{aligned}{}^t R_\hbar(y,z)\sqrt{\det g(z)} &= \hbar^{-d} b_0\left(y, \frac{y-z}{\hbar}\right) + q_\hbar(y,z), \\ \int b_0(y,u)\,du &= 1, \qquad |q_\hbar(y,z)| \le \hbar^{-d+1}\psi\left(\frac{|y-z|}{\hbar}\right).\end{aligned} \tag{2.46}$$

Set $y = h\hat{y}, z = h\hat{z}$ and $\hat{\Theta}_{h,s}(0, \hat{y}) = h^d \Theta_{h,s}(x_0, y)$. Then (2.44) becomes

$$(2.47) \quad \sup_{h \in ]0,h_0]} \int_{|\hat{y}| \leq 2c_0 h^{-1}} |\hat{\Theta}_{h,s}(0, \hat{y})| \sqrt{\det g(x_0 + h\hat{y})} \, d\hat{y} \leq C_0/\sqrt{s}.$$

One has by (2.46),

$$(2.48) \quad \left| \int q_\hbar(y, z) f_{x_0}(z) \, d_g z \right| \leq C \int \hbar^{-d+1} \psi\left(\frac{|y-z|}{\hbar}\right) \frac{\mathbf{1}_{\{d_g(x_0, z) \leq h\}}}{|B(x_0, h)|} \, dz$$
$$\leq C\hbar^{-d+1} \int_{d_g(0, h\hat{z}) \leq h} \psi(\sqrt{s}|\hat{y} - \hat{z}|) \, d\hat{z}$$

and

$$(2.49) \quad f_{x_0}(y) - \int \hbar^{-d} b_0\left(y, \frac{y-z}{\hbar}\right) f_{x_0}(z) \, dz$$
$$= \int \hbar^{-d} b_0\left(y, \frac{y-z}{\hbar}\right) (f_{x_0}(y) - f_{x_0}(z)) \, dz.$$

From (2.48) and (2.49), we get for some $\psi$ with rapid decay,

$$(2.50) \quad |\hat{\Theta}_{h,s}(0, \hat{y})| \leq C \int s^{d/2} \psi(\sqrt{s}|\hat{y} - \hat{z}|)$$
$$\times (\hbar \mathbf{1}_{\{d_g(0, h\hat{z}) \leq h\}} + |\mathbf{1}_{\{d_g(0, h\hat{z}) \leq h\}} - \mathbf{1}_{\{d_g(0, h\hat{y}) \leq h\}}|) \, d\hat{z}.$$

This implies

$$(2.51) \quad \int_{|\hat{y}| \leq 2c_0 h^{-1}} |\hat{\Theta}_{h,s}(0, \hat{y})| \sqrt{\det g(x_0 + h\hat{y})} \, d\hat{y}$$
$$\leq C \int \int s^{d/2} \psi(\sqrt{s}|\hat{y} - \hat{z}|)$$
$$\times (\hbar \mathbf{1}_{\{d_g(0, h\hat{z}) \leq h\}} + |\mathbf{1}_{\{d_g(0, h\hat{z}) \leq h\}} - \mathbf{1}_{\{d_g(0, h\hat{y}) \leq h\}}|) \, d\hat{z} \, d\hat{y}$$
$$\leq C\left(\hbar + \int_0^1 \int_{u\sqrt{s}}^\infty \frac{dv}{1+v^4} \, du\right) \leq C_0/\sqrt{s}.$$

The proof of our lemma is complete. $\square$

### 3. The spectral theory of $T_h$.

3.1. *Estimates on eigenfunctions.* In this section, we prove estimates on the eigenfunctions of $T_h$. Let us recall that $\|f\|_{H^s(M)}$ denotes the usual Sobolev norm, and that the semi-classical Sobolev norm $\|f\|_{h,s}$ is defined by (2.6). For a family $f_h \in L^2(M)$, we shall write $f_h \in \mathcal{O}_{C^\infty}(h^\infty)$ iff there exists $h_0 > 0$, such that for any $s, N$ there exists $C_{s,N}$ such that one has $\|f_h\|_{H^s(M)} \leq C_{s,N} h^N$ for all $h \in ]0, h_0]$. If $f_h = \sum f_{j,h} e_j$ is the Fourier expansion of $f_h$ in the basis of eigenfunctions of $\Delta_g$, this is equivalent to

$$\exists h_0 > 0, \forall k \, \forall N \, \exists C_{k,N} \qquad |f_{j,h}| \leq C_{k,N} h^N (1+\lambda_j)^{-k}$$
(3.1)
$$\forall j, \forall h \in ]0, h_0].$$

Let $0 < \delta < 1$ and $h_0 > 0$. For $h \in ]0, h_0]$, let $e^h$ be an eigenfunction of $T_h$ with $\|e^h\|_{L^2} = 1$, associated to an eigenvalue $z_h \in [\delta, 1]$, so $(T_h - z_h) e^h = 0$.

LEMMA 5. *There exists $h_0 > 0$, and for all $j \in \mathbb{N}$ there exists $C_j > 0$, such that, the following inequality holds true*

$$\sup_{h \in ]0, h_0]} \|e^h\|_{h,j} \leq C_j. \tag{3.2}$$

PROOF. We use the notation of Lemma 2 and we set $T_{h,k} = \psi_k T_h \varphi_k$. One has for $h$ small enough $T_h = \sum_k T_{h,k}$. For any given $k$, we denote by $x = (x_1, \ldots, x_d)$ local coordinates in $U_k$, and we choose a partition of unity in $\mathbb{R}^d$ of the form

$$1 = \sum_{\alpha \in \mathbb{Z}^d} \theta_\alpha, \qquad \theta_\alpha(x) = \theta\left(\frac{x - \alpha h}{h}\right), \tag{3.3}$$

with $\theta \in C_0^\infty$. Then, for any integer $m$, there exists $D_m$ independent of $h$ such that for any $u \in H^m(\mathbb{R}^d)$ with compact support, one has

$$D_m^{-1} \sum_\alpha \|\theta_\alpha u\|_{h,m}^2 \leq \|u\|_{h,m}^2 \leq D_m \sum_\alpha \|\theta_\alpha u\|_{h,m}^2. \tag{3.4}$$

If $\theta' \in C_0^\infty$ is equal to 1 on the set $\{X, \text{dist}(X, \text{support}(\theta)) \leq 2\}$, one has for $h \in ]0, h_0]$ with $h_0 > 0$ small enough, $\theta_\alpha T_h = \theta_\alpha T_h \theta'_\alpha$ for all $\alpha$. For any given $\alpha$, we perform the change of variable $x = h(\alpha + X)$. Let $S_\alpha$ be the rescaled operator acting on functions of the variable $X$ defined by [with $f(x) = F(\frac{x-\alpha h}{h})$]

$$\theta_\alpha T_{h,k} \theta'_\alpha (f)(h(\alpha + X)) = S_\alpha(F)(X). \tag{3.5}$$

Let us first show that $S_\alpha$ is the sum of two quantized canonical transformations of degree $-(1+d)/2 \leq -1$. From the definition (3.5) of $S_\alpha$ and (2.16), one has

$$S_\alpha(F)(X) = (2\pi)^{-d} \frac{h^d \theta(X)}{|B(h(\alpha + X), h)|}$$
(3.6)
$$\times \int e^{i(X-Y)\xi} q(h(\alpha + X), \xi, h) \theta'(Y) F(Y) \, dY \, d\xi,$$

$$q(h(\alpha + X), \xi, h) = \int_{|v| \leq 1} e^{i\xi \cdot m(h(\alpha+X), hv)v} \rho(h(\alpha+X), hv) \, dv.$$

Let us compute the integral which defined $q(x,\xi,h)$ for $|\xi|$ large. The phase $v \to \xi.m(x,hv)v$ as no critical points in $v$, so if $\chi(r) \in C_0^\infty ]0,2[$ is equal to 1 near $r=1$, one has

$$(3.7) \quad q(x,\xi,h) = \int_0^1 \chi(r) r^{d-1} \left( \int_{|\omega|=1} e^{i\xi.rm(x,hr\omega)\omega} \rho((x,hr\omega)) \, d\omega \right) dr + n(x,\xi,h),$$

where $n$ is a symbol in $S^{-\infty}$. The phase $\omega \to \xi.rm(x,hr\omega)\omega$ has two non-degenerate critical points on the sphere $|\omega|=1$, $\omega_c^\pm = \pm \frac{g^{-1/2}(x)\xi}{|g^{-1/2}(x)\xi|} + \mathcal{O}(h)$, since $\pm \frac{g^{-1/2}(x)\xi}{|g^{-1/2}(x)\xi|}$ are the two nondegenerate critical points of the phase $\omega \to \xi.rm(x,0)\omega$, and the critical values (homogeneous in $\xi$ of degree 1) are $\Phi_\pm(x,r,\xi,h) = \pm r|\xi|_x + \mathcal{O}(h)$ since $|g^{-1/2}(x)\xi| = |\xi|_x$. Using the stationary phase theorem, we get

$$(3.8) \quad q(x,\xi,h) = \int_0^1 \chi(r) r^{d-1} (e^{i\Phi_+(x,r,\xi,h)} \sigma_+(x,r,\xi,h) + e^{i\Phi_-(x,r,\xi,h)} \sigma_-(x,r,\xi,h)) \, dr + n(x,\xi,h),$$

where $\sigma_\pm$ are two symbols of degree $-(d-1)/2$. By integration in $r$, we thus get

$$(3.9) \quad q(x,\xi,h) = e^{i\Phi_+(x,1,\xi,h)} \tau_+(x,\xi,h) + e^{i\Phi_-(x,1,\xi,h)} \tau_-(x,\xi,h) + n(x,\xi,h),$$

where $\tau_\pm$ are two symbols of degree $-(d+1)/2$. From (3.9) and (3.6), we get that $S_\alpha$ is (uniformly in $\alpha, h$ for $h \in ]0,h_0]$ with $h_0 > 0$ small), the sum of two quantized canonical transformations of degree $-(d+1)/2$, with canonical relations closed to the ones associated to the phases $(X-Y)\xi \pm |\xi|_{h(\alpha+X)}$, that is, of the form $(Y,\eta) \mapsto (X = Y \pm \eta/|\eta|_{h\alpha} + \mathcal{O}(h), \xi = \eta + \mathcal{O}(h))$.

Since $T_h$ is (in the variable $X$) the sum of two quantized canonical transformations of degree $-(1+d)/2$, and since $e^h = \frac{1}{z_h} T_h(e^h)$, and $z_h \geq \delta$, we get that there exists $c$ and for all $m$, $C_m$, independent of $h, \alpha$, such that

$$(3.10) \quad \|\theta(X) e^h(h(\alpha + X))\|_{H_X^m} \leq C_m \|\theta'(X) e^h(h(\alpha + X))\|_{H_X^{m-1}},$$

where $H_X^m$ denotes the Sobolev space in variable $X$, as soon as $\theta'(X)$ is equal to 1 at each point $X$ whose distance to support($\theta$) is less than $c$. From (3.10) with $m=1$, (3.4), and $h\partial_x = \partial_X$, we get for $\chi(x) \in C_0^\infty(U_k)$ and $h \in ]0,h_0]$ with $h_0 > 0$ small,

$$(3.11) \quad \|\chi(x) e^h(x)\|_{h,1} \leq C \|e^h(x)\|_{L^2(U_k)}.$$

Therefore, since $(U_k)$ is a covering of $M$, we get $\|e^h\|_{h,1} \leq C\|e^h\|_{L^2}$. We can now iterate this argument from (3.10), and we get for any $j$,

$$\|e^h\|_{h,j} \leq C_j\|e^h\|_{L^2}. \tag{3.12}$$

The proof of our lemma is complete. $\square$

Remark that there exists $s_1 > 1$ such that $|\Gamma_d(s)| \leq \frac{\delta}{2}$ for all $s \geq s_1 - 1$. Let $\chi \in C_0^\infty(\mathbb{R}_+)$ be equal to 1 on $[0, s_1]$ and equal to 0 on $[s_1 + 1, \infty[$.

LEMMA 6. *Let $e^h$ as in Lemma 5. Then*

$$\chi(-h^2\Delta_g)e^h - e^h = \mathcal{O}_{C^\infty}(h^\infty). \tag{3.13}$$

PROOF. Let $(e_j)_{j\in\mathbb{N}}$ be an Hilbertian basis of $L^2(M, d_gx)$ such that $-\Delta_g e_j = \lambda_j e_j$ and consider $\Pi_s$ the orthogonal projector on $\text{span}\{e_j, h^2\lambda_j \geq s\}$. By Lemma 4, there exist $s_0, h_0$ such that

$$\forall s \geq s_0 \quad \sup_{h\in]0,h_0]} \|\Pi_s T_h \Pi_s\|_{L^2} \leq \delta/2. \tag{3.14}$$

Let $s_2 > \max(s_1 + 1, s_0)$ and let $\chi_2, \chi_3$ be smooth functions such that $1_{\mathbb{R}_+} = \chi + \chi_2 + \chi_3$, $\chi_3(s) = 0$ for $s \leq s_2 - 1$ and $\chi_3(s) = 1$ for $s \geq s_2$. Let $\tilde\chi_2 \in C_0^\infty(\mathbb{R}_+)$ equal to 1 near $[s_1, s_2]$ and equal to 0 on $[0, s_1 - 1] \cup [s_2 + 1, \infty[$. On $\text{supp}(\tilde\chi_2(s))$ we have $z_h - \Gamma_d(s) \geq \frac{\delta}{2}$. Hence it follows from Lemma 3 that there exist $E \in \mathcal{E}^0$ such that $E(T_h - z_h) = \tilde\chi_2(-h^2\Delta_g) + R$ with $R \in h^\infty \mathcal{E}^{-\infty}$. As $(T_h - z_h)e^h = 0$, we get

$$\tilde\chi_2(-h^2\Delta_g)e^h \in \mathcal{O}_{C^\infty}(h^\infty). \tag{3.15}$$

Set $e^h = \sum_j x_j^h e_j$. Then

$$\begin{aligned}\Pi_{s_2}e^h - \chi_3(-h^2\Delta_g)e^h &= \sum_{h^2\lambda_j \geq s_2} x_j^h e_j - \sum_j \chi_3(h^2\lambda_j)x_j^h e_j \\ &= -\sum_{s_1 \leq h^2\lambda_j < s_2} \chi_3(h^2\lambda_j)x_j^h e_j.\end{aligned} \tag{3.16}$$

As $\tilde\chi_2 = 1$ on $[s_1, s_2]$, it follows from (3.15) and (3.1) that one has $\Pi_{s_2}e^h - \chi_3(-h^2\Delta_g)e^h \in \mathcal{O}_{C^\infty}(h^\infty)$. Therefore we get

$$e^h = \chi(-h^2\Delta_g)e^h + \Pi_{s_2}e^h + \mathcal{O}_{C^\infty}(h^\infty). \tag{3.17}$$

Since $\Pi_{s_2}$ is bounded by 1 on $L^2$, applying $\Pi_{s_2}(T_h - z_h)$ to this equality, we get

$$\Pi_{s_2}(T_h - z_h)\Pi_{s_2}e^h = -\Pi_{s_2}(T_h - z_h)\chi(-h^2\Delta_g)e^h + \mathcal{O}_{L^2}(h^\infty). \tag{3.18}$$

Let $\tilde\chi \in C_0^\infty([0,\infty[)$ be supported in $[0, s_2[$ and equal to 1 near the support of $\chi$. Then, thanks to Lemma 2, we have

(3.19) $\quad (T_h - z_h)\chi(-h^2\Delta_g) = \tilde\chi(-h^2\Delta_g)(T_h - z_h)\chi(-h^2\Delta_g) + h^\infty \mathcal{E}_{cl}^{-\infty}.$

From (3.18), (3.19) and $\Pi_{s_2}\tilde\chi(-h^2\Delta_g) = 0$, we get $\Pi_{s_2}(T_h - z_h)\Pi_{s_2}e^h \in \mathcal{O}_{L^2}(h^\infty)$. Since $s_2 > s_0$, the operator $\Pi_{s_2}(T_h - z_h)\Pi_{s_2}$ is invertible on the space $\Pi_{s_2}(L^2(M))$. Consequently, $\Pi_{s_2}e^h$ is $O(h^\infty)$ in $L^2(M)$. On the other hand, from Lemma 5, (3.2), one has for any integer $j$, $\|\Delta_g^{j/2}\Pi_{s_2}e^h\|_{L^2} = \|\Pi_{s_2}\Delta_g^{j/2}e^h\|_{L^2} \leq C_j h^{-j}$. By interpolation it follows that for all $j$, one has $\|\Delta_g^{j/2}\Pi_{s_2}e^h\|_{L^2} \in O(h^\infty)$, that is, one has $\Pi_{s_2}e^h \in \mathcal{O}_{C^\infty}(h^\infty)$. Then (3.13) follows from (3.17). The proof of our lemma is complete. $\square$

For $z_h \in [\delta, 1]$, set $z_h = 1 - h^2\tau_h$, so that $e^h$ satisfies $T_h e^h = (1 - h^2\tau_h)e^h$. The next lemma is a refinement of Lemma 5.

LEMMA 7. *For all $j \in \mathbb{N}$, there exists $C_j$ such that for all $h \in ]0, h_0]$, the following inequality holds true:*

(3.20) $$\|e^h\|_{H^j(M)} \leq C_j(1 + \tau_h)^{j/2}.$$

PROOF. By Lemma 6, we have $e^h - \chi(-h^2\Delta_g)e^h \in \mathcal{O}_{C^\infty}(h^\infty)$, and therefore using also Lemma 1, we get $((T_h - 1)\chi(-h^2\Delta_g) + h^2\tau_h)e^h \in \mathcal{O}_{C^\infty}(h^\infty)$ and it follows from Lemma 3 and $(\Gamma_d - 1)(1 - \chi)(-h^2\Delta_g)e^h \in \mathcal{O}_{C^\infty}(h^\infty)$ that

(3.21) $\quad ((\Gamma_d - 1)(-h^2\Delta_g) + h^2 A_h + h^2\tau_h)e^h \in \mathcal{O}_{C^\infty}(h^\infty)$

with $A_h \in \mathcal{E}_{cl}^{-\infty}$. One has $(\Gamma_d - 1)(s) = -sF_d(s)$ with $F_d$ smooth, and from (3.21), we get

(3.22) $\quad -\Delta_g F_d(-h^2\Delta_g)e^h = (A_h + \tau_h)e^h + \mathcal{O}_{C^\infty}(h^\infty).$

Since $A_h$ is uniformly in $h$ bounded on all $H^j(M)$, and $\|e^h\|_{L^2} = 1$, we get from (3.22) for all $j \in \mathbb{N}$, with $C_j$ independent of $h$,

(3.23) $\quad \|F_d(-h^2\Delta_g)e^h\|_{H^{j+2}(M)} \leq C_j(1 + \tau_h)\|e^h\|_{H^j(M)}.$

Since $F_d(s) \neq 0$ on $[0, s_1 + 2]$, we get (3.20) by induction on $j$ from (3.23) and (3.13). The proof of our lemma is complete. $\square$

3.2. *Proof of Theorem 1.* Let us recall that there exists $\gamma_0 < 1$ such that $\Gamma_d(s) \in [-\gamma_0, 1]$ for all $s \in \mathbb{R}$. Let $\varepsilon \in ]0, (1 - \gamma_0)/2[$ and $\chi(t) \in C_0^\infty([0,\infty[)$ equal to 1 near $t = 0$ and such that $\chi(t) \in [0, 1]$ for all $t$. Thanks to Lemma 4, there exists $s > 0$ such that

(3.24) $\quad \left\|T_h(1 - \chi)\left(\dfrac{-h^2\Delta_g}{s}\right)\right\|_{L^2(M, d\nu_h)} \leq \varepsilon.$

On the other hand, thanks to Lemma 3 we can apply the Garding inequality to the pseudo-differential operator $T_h \chi(\frac{-h^2 \Delta_g}{s})$ to get for $h > 0$ small enough,

$$(3.25) \quad \left\langle T_h \chi\left(\frac{-h^2 \Delta_g}{s}\right) f, f \right\rangle_{L^2(M, d\nu_h)} \geq (-\gamma_0 - \varepsilon)\|f\|_{L^2(M, d\nu_h)},$$

where we have used the fact that $\sup_{f \neq 0} \|f\|_{L^2}/\|f\|_{L^2(M, d\nu_h)}$ goes to 1 when $h$ goes to 0. Combining equations (3.24) and (3.25), we obtain

$$(3.26) \qquad \langle T_h f, f \rangle_{L^2(M, d\nu_h)} \geq (-\gamma_0 - 2\varepsilon)\|f\|^2_{L^2(M, d\nu_h)}$$

which proves the first statement of Theorem 1 as $T_h$ is self-adjoint on $L^2(M, d\nu_h)$.

Let us now prove (1.8). Set $|\Delta_h| = 2(d+2)\frac{1-T_h}{h^2}$. For $k \leq L$, we denote by $m_k = \dim(Ker(\Delta_g + \lambda_k))$ the multiplicity of $\lambda_k$. Let $\rho_0 \in C_0^\infty(\mathbb{R})$ be equal to 1 near zero. Then there exists $h_0 > 0$ such that for $h \in \,]0, h_0]$, one has $e = \rho_0(-h^2 \Delta_g)e$ for any $e \in Ker(\Delta_g + \lambda_k)$ with $k \leq L$. Thus, if $(U_j)$ is a finite covering of $M$ by local charts and $1 = \sum \varphi_j$ a partition of unity with $\varphi_j \in C_0^\infty(U_j)$, one has

$$(3.27) \qquad (T_h - \Gamma_{d,h})(e) = \sum_j (T_h - \Gamma_{d,h}) \rho_0(-h^2 \Delta_g) \varphi_j(e).$$

From Lemma 3 one has for each $j$, $(T_h - \Gamma_{d,h})\rho_0(-h^2 \Delta_g)\varphi_j = h^2 Op(a) + R_h$, with $a = a_2 + ha_3 + \cdots \in S_{cl}^{-\infty}$ compactly supported in $x \in U_j$, $R_h$ smoothing, $a_2(x, \xi) = \mathcal{O}(\xi^2)$ near $\xi = 0$ and $a_3(x, 0) = 0$. As $e$ is smooth and does not depend on $h$, it follows that $((T_h - \Gamma_{d,h})\rho_0(-h^2 \Delta_g)\varphi_j(e) \in O_{L^2}(h^4)$. Therefore,

$$(3.28) \qquad \|(T_h - \Gamma_{d,h})(e)\|_{L^2(M, d\nu_h)} = O(h^4).$$

Moreover, $\Gamma_{d,h} e = \Gamma_d(h^2 \lambda_k) e = (1 + h^2 \Gamma'_d(0)\lambda_k + O(h^4))e$. Combining this with (3.28) we obtain $\|(|\Delta_h| - \lambda_k)e\|_{L^2(M, d\nu_h)} = O(h^2)$ for all $e \in Ker(\Delta_g + \lambda_k)$, and since $|\Delta_h|$ is self-adjoint on $L^2(M, d\nu_h)$, we get that there exists $C_0$ such that

$$(3.29) \quad \begin{aligned} &\forall h \in \,]0, h_0], \forall 0 \leq k \leq L \\ &\operatorname{card}(Spec(|\Delta_h|) \cap [\lambda_k - C_0 h^2, \lambda_k + C_0 h^2]) \geq m_k. \end{aligned}$$

Now, if $e^h$ is a normalized eigenfunction of $|\Delta_h|$, $|\Delta_h| e^h = \tau_h e^h$, with $\tau_h$ bounded, one has, by Lemma 6, $e^h - \rho_0(-h^2 \Delta_g)e^h \in \mathcal{O}_{C^\infty}(h^\infty)$, and also by Lemma 7 since $\tau_h$ is bounded, $\|e^h\|_{H^j(M)} \leq C_j$ for all $j$, with $C_j$ independent of $h$. Thus the same argument as above shows that there exists $C$ independent of $h$ such that

$$(3.30) \qquad \|(\tau_h + \Delta_g)(e^h)\|_{L^2(M, d\nu_h)} \leq C h^2,$$

and thus $\text{dist}(\tau_h, Spec(-\Delta_g)) \leq Ch^2$. It remains to prove that for $h$ small, we have equality in the right-hand side of (3.29). Let $p \geq m_k$ and let $e_1(h), \ldots, e_p(h)$ be a family of eigenfunctions of $|\Delta_h|$ associated to the eigenvalues $\tau_j(h) \in [\lambda_k - C_0 h^2, \lambda_k + C_0 h^2]$, orthonormal for the scalar product $\langle \cdot, \cdot \rangle_{L^2(M, d\nu_h)}$. By Lemma 7, there exists a sequence $(h_n)$ going to zero as $n \to \infty$ such that $e_l(h_n)$ converges in $H^2$. Denoting $f_l$ its limit we get from (3.30), $-\Delta_g f_l = \lambda_k f_l$ for all $l = 1, \ldots, p$ and the functions $f_l$ are orthogonal for the scalar product $\langle \cdot, \cdot \rangle_{L^2(M, d_g x)}$. This proves that $m_k \geq p$, and completes the proof of (1.8). (In particular, this implies that 1 is a simple eigenvalue of $T_h$.)

Let us now prove the Weyl estimate (1.9).

Let $\delta \in ]0, 1[$ be given. Let $\tau \in [0, (1-\delta)h^{-2}]$. Observe that $N(1 - \tau h^2, h)$ is the number of eigenvalues of $|\Delta_h|$ in the interval $[0, 2(d+2)\tau]$. We denote by $N_0(a, h)$ the number of eigenvalues of $\Gamma_d(-h^2 \Delta_g)$ in the interval $[a, 1]$. Let us define the function $\Phi_h(s)$ and the operator $|\Delta_h^0|$ by the formulas,

$$\Phi_h(s) = 2(d+2)\frac{1 - \Gamma_d(s)}{h^2},$$
(3.31)
$$|\Delta_h^0| = \Phi_h(-h^2 \Delta_g).$$

Then $N_0(1 - \tau h^2, h)$ is the number of eigenvalues of $|\Delta_h^0|$ in the interval $[0, 2(d+2)\tau]$. Let us first show that $N_0$ satisfies the Weyl estimate (1.9), that is, there exists $C$ such that for all $h \in ]0, h_0]$ and all $\tau \in [0, (1-\delta)h^{-2}]$, one has

$$(3.32) \quad \left| N_0(1 - \tau h^2, h) - (2\pi h)^{-d} \int_{\Gamma_d(|\xi|_x^2) \in [1 - \tau h^2, 1]} dx\, d\xi \right|$$
$$\leq C(1 + \tau)^{(d-1)/2}.$$

To prove this point, let $n^+(\lambda)$ [resp. $n^-(\lambda)$] be the number of eigenvalues $\lambda_j$ of $-\Delta_g$ in the interval $[0, \lambda]$ (resp. $[0, \lambda[$). By the classical Weyl estimate with accurate remainder (see [7]), one has

$$(3.33) \quad n^{\pm}(\lambda) = (2\pi)^{-d} \int_{|\xi|_x^2 \leq \lambda} dx\, d\xi + \mathcal{O}(\lambda^{(d-1)/2}).$$

By (3.31), $N_0(1 - \tau h^2, h)$ is the number of eigenvalues $\lambda_j$ of $-\Delta_g$ such that $1 - \Gamma_d(h^2 \lambda_j) \leq \tau h^2$. Since $\tau \leq (1-\delta)h^{-2}$, the set $\{s \geq 0; 1 - \Gamma_d(s) \leq \tau h^2\}$ is a finite union of disjoint intervals $I_0 \cup \cdots \cup I_k$ with $I_0 = [0, s_0(\tau h^2)]$, $I_j = [s_j^-(\tau h^2), s_j^+(\tau h^2)]$ for $1 \leq j \leq k$, and such that $c_0 \leq s_1^- \leq s_1^+ < s_2^- \leq s_2^+ < \cdots \leq s_k^+ \leq c_1$ with $c_0 > 0$ independent of $h, \delta$ and $c_1$ independent of $h$. Thus we get

$$(3.34) \quad N_0(1 - \tau h^2, h) = n^+(s_0 h^{-2}) + \sum_{j=1}^{j=k} n^+(s_j^+ h^{-2}) - n^-(s_j^- h^{-2}).$$

Observe that $k = 0$ when $\tau \leq ch^{-2}$ with $c$ small enough, and in that case one has by (1.6), $s_0 h^{-2} \simeq 2(d+2)\tau$, and therefore (3.32) is consequence of (3.33). On the other hand, in the case $\tau \geq ch^{-2}$, then both $(s_j^\pm h^{-2})^{(d-1)/2}$ and $(s_0 h^{-2})^{(d-1)/2}$ are of order $\tau^{(d-1)/2}$, and thus we get (3.32) from (3.33) and (3.34).

Let $E_\tau$ be the finite dimension space spanned by the eigenfunctions $e_j$ of $-\Delta_g$ with $\Phi_h(h^2\lambda_j) \leq 2\tau(d+2)$. Then by (3.31), one has $\dim(E_\tau) = N_0(1 - \tau h^2, h)$. By (2.30) and $\||\Delta_h|\|_{L^2} \leq Ch^{-2}$, one has for all $f \in L^2$,

$$(3.35) \quad |(|\Delta_h|f|f)_{L^2(M, Z_h d\nu_h)} - (|\Delta_h|f|f)_{L^2(M, d_g x)}| \leq C' \|f\|_{L^2}^2.$$

Let $\chi \in C_0^\infty([0, \infty[)$ equal to 1 near the compact set $\{s \geq 0; 1 - \Gamma_d(s) \leq 1 - \delta\}$. Then $f = \chi(-h^2\Delta_g)f$ for all $f \in E_\tau$, and from Lemma 3, one has $(|\Delta_h| - |\Delta_h^0|)\chi(-h^2\Delta_g) = -2(d+2)A_h$. Thus, since $A_h$ is bounded on $L^2$, from (3.35) we get that there exists $C_- = C_-(\delta)$ independent of $\tau, h$, such that for all $f \in E_\tau$, one has

$$(3.36) \quad (|\Delta_h|f|f)_{L^2(M, Z_h d\nu_h)} \leq 2(\tau + C_-)(d+2)\|f\|_{L^2(M, Z_h d\nu_h)}^2,$$

and this implies, by the min-max,

$$(3.37) \quad N_0(1 - \tau h^2, h) = \dim(E_\tau) \leq N(1 - (\tau + C_-)h^2, h).$$

Let $F_\tau$ be the orthogonal complement of $E_\tau$ in $L^2(M, d_g x)$. Let $\theta \in C_0^\infty$ such that $\|T_h(1-\theta)(-h^2\Delta_g)\|_{L^2} \leq \delta$. Let $\chi \in C_0^\infty$ with values in $[0, 1]$, equal to 1 near $[0, 1-\delta] \cup \text{support}(\theta)$. Let $\psi = 1 - \chi$, so that $(1-\theta)\psi = \psi$. Let $A_h = (|\Delta_h| - |\Delta_h^0|)\chi(-h^2\Delta_g) \in \mathcal{E}_{cl}^{-\infty}$ and $B_h = \chi(-h^2\Delta_g)(|\Delta_h| - |\Delta_h^0|) \in \mathcal{E}_{cl}^{-\infty}$ be given by Lemma 3. Then, one has

$$(3.38) \quad |\Delta_h| = \chi|\Delta_h^0|\chi + \psi|\Delta_h^0|\chi + \chi|\Delta_h^0|\psi + \psi|\Delta_h|\psi + A_h + B_h\psi.$$

The operator $A_h + B_h\psi$ is bounded on $L^2$ by a constant $C(\delta)$ uniformly in $h$. From $\psi(1 - T_h)\psi = \psi^2 - \psi T_h(1-\theta)\psi$, we get

$$(3.39) \quad (\psi|\Delta_h|\psi f|f)_{L^2(M, d_g x)} \geq 2(1-\delta)\frac{d+2}{h^2}\|\psi f\|_{L^2(M, d_g x)}^2.$$

Therefore, from (3.35) we get that there exists $C_+ = C_+(\delta) > 0$ independent of $\tau, h$, such that for all $f = \sum_{\lambda_j > \tau} x_j e_j \in F_\tau$, one has

$$(3.40) \quad \begin{aligned} &(|\Delta_h|f|f)_{L^2(M, Z_h d\nu_h)} + (d+2)C_+\|f\|_{L^2(M, Z_h d\nu_h)}^2 \\ &\geq \sum_{\lambda_j > \tau} \Phi_h(h^2\lambda_j)(\chi^2 + 2\chi\psi)(h^2\lambda_j)|x_j|^2 \\ &\quad + 2(1-\delta)\frac{d+2}{h^2}\sum_{\lambda_j > \tau} \psi^2(h^2\lambda_j)|x_j|^2 \\ &\geq 2\tau(d+2)\sum_{\lambda_j > \tau} |x_j|^2 \geq (2\tau(d+2) - Ch^2)\|f\|_{L^2(M, Z_h d\nu_h)}^2, \end{aligned}$$

and this implies by the min–max for $\tau$ large enough, and $h \in ]0, h_0]$ with $h_0$ small,

$$(3.41) \qquad N_0(1 - \tau h^2, h) = \operatorname{codim}(F_\tau) \geq N(1 - (\tau - C_+)h^2, h).$$

Then we obtain the Weyl estimate (1.9) from (3.32), (3.37) and (3.41). Finally, the estimate (1.11) is an easy byproduct of the estimates (3.20) of Lemma 7. The proof of Theorem 1 is complete.

3.3. *Proof of Theorem 2.* Let us recall that $\Phi_h(s)$ and $|\Delta_h^0|$ are defined in (3.31).

One has $2(d+2)(1 - \Gamma_d(s)) \geq c_1 \min(s, 1)$ with $c_1 > 0$, and, therefore,

$$(3.42) \qquad \Phi_h(h^2 \lambda_j) \geq c_1 \min(\lambda_j, h^{-2}).$$

Observe that there exists $h_0, c_0 > 0$ such that for all $z \in U$, all $h \in ]0, h_0]$, and all $j \in \mathbb{N}$, one has

$$(3.43) \qquad |z - \Phi_h(h^2 \lambda_j)| \geq c_0(1 + |z| + \min(\lambda_j, h^{-2})).$$

To see this fact, just observe that by (3.42), for $c_1 \min(\lambda_j, h^{-2}) \geq A + 1$, (3.43) holds true, since $z \in U$. Now, $c_1 \min(\lambda_j, h^{-2}) \leq A + 1$ implies if $h_0$ is small, $\lambda_j \leq (A+1)/c_1$, and therefore, $|\Phi_h(h^2 \lambda_j) - \lambda_j| \leq c_2 h^2$, and (3.43) holds true also in that case since $z \in U$. Since for $h^2 \lambda_j \leq c_3$ with $c_3 > 0$ small, one has $|\Phi_h(h^2 \lambda_j) - \lambda_j| \leq c_4 h^2 \lambda_j^2$, we get from (3.43), that there exists $C$ such that for all $z \in U$ and all $h \in ]0, h_0]$, one has

$$(3.44) \qquad \sup_{j \in \mathbb{N}} \left| \frac{1}{z - \Phi_h(h^2 \lambda_j)} - \frac{1}{z - \lambda_j} \right| \leq C h^2,$$

and this implies, obviously,

$$(3.45) \qquad \|(z - |\Delta_h^0|)^{-1} - (z + \Delta_g)^{-1}\|_{L^2} \leq C h^2,$$

and thus we are reduced to prove the estimate

$$(3.46) \qquad \|(z - |\Delta_h|)^{-1} - (z - |\Delta_h^0|)^{-1}\|_{L^2} \leq C h^2.$$

Observe that, as a straightforward consequence of Theorem 1 and of the self-adjointness of $|\Delta_h|$ and $|\Delta_h^0|$, respectively, on $L^2(M, d\nu_h)$ and $L^2(M, d_g x)$, there exists $C > 0$ and $h_0 > 0$ such that for all $z \in U$ and all $h \in ]0, h_0]$,

$$(3.47) \qquad \|(z - |\Delta_h|)^{-1}\|_{L^2} + \|(z - |\Delta_h^0|)^{-1}\|_{L^2} \leq \frac{C}{1 + |z|}.$$

Therefore, in order to prove (3.46), we may, and will assume that $z$ satisfies $h^2 |z| \leq \alpha$, with $\alpha > 0$ small. Using Lemma 4, we then choose $\chi_0 \in C_0^\infty$ equal to 1 on $[0, s_0]$, with support in $[0, 2s_0]$, and, such that,

$$(3.48) \qquad \|2(d+2) T_h (1 - \chi_0)(-h^2 \Delta_g)\|_{L^2} \leq d + 2 - \alpha/2.$$

Let $\chi \in C_0^\infty$ equal to 1 near $[0, 3s_0]$, and set $R_h = (z - |\Delta|_h)^{-1} - (z - |\Delta|_h^0)^{-1}$. Then since $|\Delta_h^0|$ commutes with $\Delta_g$, one has

$$
(3.49) \quad \begin{aligned} R_h \chi(-h^2 \Delta_g) \\ = (z - |\Delta_h|)^{-1}(|\Delta_h| - |\Delta_h^0|)\chi(-h^2 \Delta_g)(z - |\Delta_h^0|)^{-1}. \end{aligned}
$$

From Lemma 3, one has $(|\Delta_h| - |\Delta_h^0|)\chi(-h^2 \Delta_g) = A_h \chi'(-h^2 \Delta_g)$, with $\chi'$ equal to 1 near the support of $\chi$, and the operator $A_h \in \mathcal{E}_{cl}^{-\infty}$ satisfies

$$(3.50) \quad \|A_h f\|_{L^2(M)} \leq C h^2 \|f\|_{H^2(M)}.$$

On the other hand, from (3.43), we get

$$(3.51) \quad \|\chi'(-h^2 \Delta_g)(z - |\Delta_h^0|)^{-1} f\|_{H^2(M)} \leq C \|f\|_{L^2(M)}.$$

From (3.47), (3.49), (3.50) and (3.51), we get

$$(3.52) \quad \|R_h \chi(-h^2 \Delta_g)\|_{L^2} \leq C h^2.$$

It remains to estimate $R_h(1 - \chi)(-h^2 \Delta_g)$, and it is obviously sufficient to prove the two estimates

$$(3.53) \quad \|(z - |\Delta_h|)^{-1}(1 - \chi)(-h^2 \Delta_g)\|_{L^2} \leq C h^2,$$

$$(3.54) \quad \|(z - |\Delta_h^0|)^{-1}(1 - \chi)(-h^2 \Delta_g)\|_{L^2} \leq C h^2.$$

Since $\chi(s) = 1$ near $s = 0$, (3.54) is a consequence of (3.43). Let $g \in L^2(M)$ with $\|g\|_{L^2} = 1$ and let $f = (z - |\Delta_h|)^{-1}(1 - \chi)(-h^2 \Delta_g) g$. Then

$$(3.55) \quad (h^2 z - 2(d+2)(1 - T_h))f = h^2 (1 - \chi)(-h^2 \Delta_g) g.$$

Let $\chi_1, \chi_2 \in C_0^\infty$ with support in $[0, 3s_0[$, with $\chi_2$ equal to 1 near the support of $\chi_1$. One has $\chi_1(1 - \chi) = 0$, and thus, multiplying (3.55) by $\chi_1(-h^2 \Delta_g)$ and using Lemma 3, we obtain

$$(3.56) \quad h^2(z - |\Delta_h^0|)\chi_1(-h^2 \Delta_g) f = h^2 A_h \chi_2(-h^2 \Delta_g) f + \mathcal{O}_{C^\infty}(h^\infty).$$

Since on the support of $\chi_1$, one has $h^2 \lambda_j \leq 3s_0$, we get from (3.43), (3.47) and (3.56) that one has $\|\chi_1(-h^2 \Delta_g)f\|_{H^2} \leq C$; thus, since $\chi_1$ is arbitrary, $\|\chi_2(-h^2 \Delta_g)f\|_{H^2} \leq C$, and from (3.56) and (3.50), we thus get

$$(3.57) \quad \|\chi_1(-h^2 \Delta_g) f\|_{L^2} \leq C h^2.$$

Then, we deduce from (3.55) and (3.57)

$$(3.58) \quad (h^2 z - 2(d+2) + 2(d+2) T_h (1 - \chi_0(-h^2 \Delta_g))) f \in \mathcal{O}_{L^2}(h^2),$$

and from (3.48), we get $\|f\|_{L^2} \leq C h^2$. The proof of Theorem 2 is complete.

## 4. Proof of Theorem 3.

4.1. *The spectral theory of the Metropolis kernel.* In this section, we will deduce from the results of Section 3, useful properties on the spectral theory of the Metropolis operator $M_h$. Let us write

$$M_h = T_h + R_h. \tag{4.1}$$

Then from (1.16) and (1.17), one has

$$\begin{aligned}R_h(f)(x) &= m_h(x)f(x) \\ &\quad + \int_{d_g(x,y)\leq h} \min\left(\frac{1}{|B(y,h)|} - \frac{1}{|B(x,h)|}, 0\right) f(y)\, d_g y.\end{aligned} \tag{4.2}$$

Let $a(x,y,h) \leq 0$ be the function

$$a(x,y,h) = h^{d-2} \min\left(\frac{1}{|B(y,h)|} - \frac{1}{|B(x,h)|}, 0\right). \tag{4.3}$$

Then $a$ is a Lipschitz function in $x$ and $y$, and from (2.30), we get that there exists $C$ independent of $x, y, h$ such that

$$|a(x,y,h)| \leq C d_g(x,y), \qquad |\nabla_x a(x,y,h)| + |\nabla_y a(x,y,h)| \leq C. \tag{4.4}$$

Since $R_h(1) = 0$, one has $m_h(x) = -h^{2-d}\int_{d_g(x,y)\leq h} a(x,y,h)\, d_g y$, and therefore the function $m_h$ is Lipschitz and satisfies $\|m_h\|_{L^\infty} \leq Ch^3$ and $\|\nabla m_h\|_{L^\infty} \leq Ch^2$. From these facts, one easily gets that there exists $C$ independent of $p \in [1,\infty]$ and $h$ such that

$$\begin{aligned}\|R_h\|_{L^p} &\leq Ch^3, \\ \|R_h\|_{W^{1,p}} &\leq Ch^2,\end{aligned} \tag{4.5}$$

where $W^{1,p} = \{f \in L^p, \nabla f \in L^p\}$ is the usual Sobolev space. Therefore, $M_h$ is a small perturbation of $T_h$. In particular, there still exist $h_0 > 0$ and $\gamma < 1$ such that the spectrum of $M_h$ is a subset of $[-\gamma, 1]$, $1$ is a simple eigenvalue of $M_h$ and since $\|m_h\|_{L^\infty} \leq Ch^3$ and $m_h(x) \geq 0$, the spectrum of $M_h$ is discrete outside $[0, Ch^3]$. Let

$$Ch^3 < \cdots \leq \widetilde{\mu}_{k+1}(h) \leq \widetilde{\mu}_k(h) \leq \cdots \leq \widetilde{\mu}_1(h) < \widetilde{\mu}_0(h) = 1 \tag{4.6}$$

be the decreasing sequence of positive eigenvalues of $M_h$. Set

$$1 - M_h = \frac{h^2}{2(d+2)}|\widetilde{\Delta}_h|. \tag{4.7}$$

Then from (4.5), one has

$$\||\widetilde{\Delta}_h| - |\Delta_h|\|_{L^2} \leq Ch. \tag{4.8}$$

From Theorem 1 and (4.8) we get that for any given $L > 0$, there exists $C$ such that for all $h \in ]0, h_0]$ and all $k \leq L$, one has

$$\left| \frac{1 - \widetilde{\mu}_k(h)}{h^2} - \frac{\lambda_k}{2(d+2)} \right| \leq Ch. \tag{4.9}$$

Moreover, since $\|T_h - M_h\|_{L^2} \leq Ch^3$, the Weyl estimate (1.9) remains valid for the number $\widetilde{N}(a, h)$ of eigenvalues of $M_h$ in the interval $[a, 1]$: for $\delta \in ]0, 1[$, one has

$$\left| \widetilde{N}(1 - \tau h^2, h) - (2\pi h)^{-d} \int_{\Gamma_d(|\xi|_x^2) \in [1-\tau h^2, 1]} dx\, d\xi \right| \tag{4.10}$$

$$\leq C_{\delta,1}(1 + \tau)^{(d-1)/2}$$

for any $\tau \in [0, (1-\delta)h^{-2}]$, and therefore, the estimate (1.10) is still valid; for any $\tau \in [0, (1-\delta)h^{-2}]$, one has

$$\widetilde{N}(1 - \tau h^2, h) \leq C_\delta(1 + \tau)^{d/2}. \tag{4.11}$$

The main result of this section is to prove that there exist $C_\delta$ such that for any eigenfunction $\widetilde{e}_k^h$ of $M_h$ associated to the eigenvalue $\widetilde{\mu}_k(h) \in [\delta, 1]$, the inequality (1.11) still holds true, that is, with $\widetilde{\tau}_k(h) = h^{-2}(1 - \widetilde{\mu}_k(h))$, one has

$$\|\widetilde{e}_k^h\|_{L^\infty} \leq C_\delta(1 + \widetilde{\tau}_k(h))^{d/4} \|\widetilde{e}_k^h\|_{L^2}. \tag{4.12}$$

We will obtain this estimate as a consequence of (4.5), using Sobolev inequalities and the following lemma.

LEMMA 8. *Let $N \geq 1$, $p \in [1, \infty]$ and $\delta \in ]0, 1[$. Let $s_0 > 0$ such that $|\Gamma_d(s)| \leq \delta/2$ for $s \geq s_0$. Let $\chi_0 \in C_0^\infty$ such that $\chi_0(s) = 1$ on $[0, s_0]$. There exist $C, C_N, h_0$, and for all $z \in K = \{z \in \mathbb{C}, |z| \in [\delta, 2]\}$ and all $h \in ]0, h_0]$, operators $E_{z,h}, \mathcal{N}_{z,h}$ which satisfy*

$$E_{z,h}(T_h - z) = 1 - \chi_0(-h^2 \Delta_g) + \mathcal{N}_{z,h}, \tag{4.13}$$

*and such that the following estimates holds true:*

$$\begin{aligned} \|E_{z,h}\|_{L^p} &\leq C, & \|E_{z,h}\|_{W^{1,p}} &\leq C, \\ \|\mathcal{N}_{z,h}\|_{L^p} &\leq C_N h^N, & \|\mathcal{N}_{z,h}\|_{W^{1,p}} &\leq C_N h^N. \end{aligned} \tag{4.14}$$

PROOF. Let $\chi \in C_0^\infty([0, 2[)$ equal to 1 on $[0, 1]$, and set $\chi_s(t) = \chi(t/s)$. By Lemma 4, there exist $s_0$ such that for all $s \geq s_0$, one has $\|T_h(1 - \chi_s(-h^2 \Delta_g))\|_{L^p} \leq \delta/2$. We then take $s \geq s_0$ such that $\chi_s = 1$ near the support of $\chi_0$, and we set $\psi = 1 - \chi_s$ and $\psi' = 1 - \chi_{4s}$. For $z \in K$, $T_h \psi - z$ is then invertible on $L^p$. Set

$$E_1 = \psi'(T_h \psi - z)^{-1}. \tag{4.15}$$

Then, there exists $C, h_0$ such that for all $h \in \,]0, h_0]$ and all $z \in K$ one has

$$\|E_1\|_{L^p} + \|E_1\|_{W^{1,p}} \leq C. \tag{4.16}$$

The $L^p$ bound is obvious since operators in $\mathcal{E}_{cl}^{-\infty}$ are bounded on $L^p$ and $\psi' = 1 - \chi_{4s}$; let us prove the $W^{1,p}$ bound in (4.16). We denote by $B$ any operator which is, uniformly in $h > 0$ small, and $z \in K$, bounded on $L^p$. Let $X$ be a vector field on $M$. Then by (2.16), one has $[T_h, X] = hB_1 X + B_2$. Thus, with $L = T_h \psi - z$, we get $[L, X] = hB_3 X + B_4$ and $[X, L^{-1}] = hB_5 X L^{-1} + B_6$. Since for $h$ small, $1 - hB_5$ is invertible on $L^p$, we obtain $XL^{-1} = B_7 X + B_8$, and thus (4.16) holds true, since $E_1 = \psi' L^{-1}$. Let $\phi \in C_0^\infty([0, 3s[)$; from $\psi' \phi = 0$, we get $E_1 L \phi = 0$, and therefore

$$E_1 \phi = E_1 [\phi, L] L^{-1}. \tag{4.17}$$

By Lemma 3, one has $[\phi, L] \in h\mathcal{E}_{cl}^{-\infty}$. Thus (4.17) implies $\|E_1 \phi\|_{L^p} + \|E_1 \phi\|_{W^{1,p}} \leq Ch$, and since $\phi$ is arbitrary, by an easy induction from (4.17), we get $\|E_1 \phi\|_{L^p} + \|E_1 \phi\|_{W^{1,p}} \leq C_N h^N$ for all $N$. Thus one has

$$E_1(T_h - z) = \psi' + \mathcal{N}_1 \tag{4.18}$$

with $\mathcal{N}_1 = E_1 T_h (1 - \psi) = E_1(\phi T_h \chi_s + \mathcal{O}(h^\infty \mathcal{E}_{cl}^{-\infty}))$ if $\phi = 1$ near $[0, 2s]$. Thus $\mathcal{N}_1$ satisfies for all $N$,

$$\|\mathcal{N}_1\|_{L^p} + \|\mathcal{N}_1\|_{W^{1,p}} \leq C_N h^N. \tag{4.19}$$

Now, by the symbolic calculus, there exist $E_2 \in \mathcal{E}_{cl}^{-\infty}$ and $\mathcal{N}_2 \in h^\infty \mathcal{E}_{cl}^{-\infty}$ such that

$$E_2(T_h - z) = \chi_{4s} - \chi_0 + \mathcal{N}_2. \tag{4.20}$$

Here we use Lemma 4 and the fact that $T_h - z$ is elliptic near the support of $\chi_{4s} - \chi_0$. Then $E_{z,h} = E_1 + E_2$ and $\mathcal{N}_{z,h} = \mathcal{N}_1 + \mathcal{N}_2$ satisfies (4.13) and (4.14). The proof of our lemma is complete. □

Let us now achieve the proof of (4.12). Let $\widetilde{\mu}(h) \in [\delta, 1]$ and $\|\widetilde{e}^h\|_{L^2} = 1$. Then $(M_h - \widetilde{\mu}(h))\widetilde{e}^h = 0$ is equivalent to $(T_h - \widetilde{\mu}(h) + R_h)\widetilde{e}^h = 0$, and using Lemma 8, we get

$$(1 - \chi_0)\widetilde{e}^h + (\mathcal{N}_{\widetilde{\mu}(h),h} + E_{\widetilde{\mu}(h),h} R_h)\widetilde{e}^h = 0. \tag{4.21}$$

Set $\widetilde{e}_l = \chi_0(\widetilde{e}^h)$ and $\widetilde{e}_+ = (1 - \chi_0)(\widetilde{e}^h)$ so that $\widetilde{e}^h = \widetilde{e}_l + \widetilde{e}_+$. Since by (4.5) and (4.13) the operator $\mathcal{N}_{\widetilde{\mu}(h),h} + E_{\widetilde{\mu}(h),h} R_h$ is $\mathcal{O}(h^2)$ on $L^p$ and $W^{1,p}$, we can solve equation (4.21) for $\widetilde{e}_+$ on the form

$$\widetilde{e}_+ = S_{\widetilde{\mu}(h),h}(\widetilde{e}_l), \tag{4.22}$$
$$\|S_{\widetilde{\mu}(h),h}\|_{L^p} + \|S_{\widetilde{\mu}(h),h}\|_{W^{1,p}} \leq Ch^2.$$

Let $1-h^2\tau = \widetilde{\mu}(h)$ and $\omega = \sqrt{1+\tau}$. One has $|\Delta_h|(\widetilde{e}^h) = 2(d+2)(\tau+h^{-2}R_h)(\widetilde{e}^h)$, and therefore, with $(|\Delta_h| - |\Delta_h^0|)\chi_0 = A_h$, we get the equation

$$(4.23) \quad |\Delta_h^0|\chi_0(\widetilde{e}^h) = (2(d+2)\chi_0(\tau+h^{-2}R_h) - A_h + [|\Delta_h|,\chi_0])(\widetilde{e}^h).$$

By (4.5) and Lemma 3, the operator $2(d+2)\chi_0(\tau+h^{-2}R_h) - A_h + [|\Delta_h|,\chi_0])$ is bounded by $\mathcal{O}(\omega^2)$ on $L^p$, uniformly in $h$. Then by (4.22) and (4.23), we get for some $p_\star \in ]d,\infty[$ and all $p \in [2,p_\star]$, that the following estimates holds true, with $C$ independent of $h$:

$$(4.24) \quad \begin{aligned} \|\widetilde{e}^h\|_{L^p} &\leq C\omega^{d/2-d/p}, \\ \|\widetilde{e}^h\|_{W^{1,p}} &\leq C\omega^{d/2-d/p+1}. \end{aligned}$$

Indeed, by (3.31) and (3.42), for $\chi_1 \in C_0^\infty$ equal to 1 near the support of $\chi_0$, one has $|\Delta_h^0|\chi_1 = -\Delta_g B_h$ with $B_h \in \mathcal{E}_{cl}^{-\infty}$ elliptic near the support of $\chi_0$. Thus, $\|\widetilde{e}^h\|_{L^2} = 1$ and (4.23) implies $\|\widetilde{e}_l\|_{W^{2,2}} \leq C\omega^2$, and thus $\|\widetilde{e}_l\|_{W^{1,2}} \leq C\omega$, so using (4.22), one gets that (4.24) holds true for $p=2$. This also shows easily that (4.24) holds true for $d=2$. When $d \geq 3$, then if (4.24) holds true for some $p \in [2,d[$, then let $q \in ]p,\infty[$ be defined by $d/q = d/p - 1$. Then the injection $W^{1,p} \subset L^q$ shows that the first line of (4.24) holds true for $q$. Moreover, in (4.23), classical properties of $-\Delta_g$ and the fact that operators in $\mathcal{E}_{cl}^{-\infty}$ are bounded on $W^{s,p}$, shows that $\|\widetilde{e}_l\|_{W^{2,p}} \leq C\omega^{d/2-d/p+2}$. Then the injection $W^{2,p} \subset W^{1,q}$ and (4.22) implies that the second line of (4.24) holds true for $q$. Then, from (4.24), we conclude the proof of (4.12) by the interpolation inequality for $p_\star > d$,

$$(4.25) \quad \|u\|_{L^\infty} \leq C\|u\|_{L^{p_\star}}^{1-d/p_\star}\|u\|_{W^{1,p_\star}}^{d/p_\star}.$$

4.2. *The total variation estimate.* In this section, we prove Theorem 3. Let $\Pi_0$ be the orthogonal projector in $L^2(M,d\mu_M)$ on the space of constant functions

$$(4.26) \quad \Pi_0(f)(x) = \frac{1}{Vol(M)}\int_M f(y)\,d_g y.$$

Then

$$(4.27) \quad 2\sup_{x \in M}\|M_h^n(x,dy) - d\mu_M\|_{\mathrm{TV}} = \|M_h^n - \Pi_0\|_{L^\infty \to L^\infty}.$$

Thus, we have to prove that there exist $A, h_0$, such that for any $n$ and any $h \in ]0,h_0]$, one has

$$(4.28) \quad e^{-\gamma'(h)nh^2} \leq \|M_h^n - \Pi_0\|_{L^\infty \to L^\infty} \leq Ae^{-\gamma(h)nh^2}$$

with $\gamma(h) \simeq \gamma'(h) \simeq \frac{\lambda_1}{2(d+2)}$ when $h \to 0$. Since $(M_h^n - \Pi_0)(\widetilde{e}_1^h) = (1-h^2\widetilde{\tau}_1^h)^n\widetilde{e}_1^h$, with $|\widetilde{\tau}_1^h - \frac{\lambda_1}{2(d+2)}| \leq Ch$ by (4.9), the lower bound in (4.28) is obvious, and to prove the upper bound, we may assume $n \geq C_0 h^{-2}$. Let $\delta \in \,]0,1[$ be such that the spectrum of $M_h$ is contained in $[-\delta, 1]$. Then write $M_h - \Pi_0 = M_{h,1} + M_{h,2}$ with

$$M_{h,1}(x,y) = \sum_{\delta \leq \widetilde{\mu}_k(h) < 1} (1 - h^2 \widetilde{\tau}_k(h)) \widetilde{e}_k^h(x) \widetilde{e}_k^h(y),$$

(4.29)
$$M_{h,2} = M_h - \Pi_0 - M_{h,1}.$$

Here $1 - h^2 \widetilde{\tau}_k(h) = \widetilde{\mu}_k(h)$. One has $M_h^n - \Pi_0 = M_{h,1}^n + M_{h,2}^n$, and we will get the upper bound in (4.28) for each of the 2 terms. From (4.29) and (4.12), there exist some $\alpha > 0$ such that

(4.30) $\quad \|M_{h,1}^n\|_{L^\infty \to L^\infty} \leq \sum_{\widetilde{\tau}_1(h) \leq \widetilde{\tau}_k(h) \leq (1-\delta)h^{-2}} (1 - h^2 \widetilde{\tau}_k(h))^n (1 + \widetilde{\tau}_k(h))^\alpha.$

Using $1 - x \leq e^{-x}$, and the estimate (4.11) on the number of eigenvalues of $M_h$ in $[1 - h^2 \tau, 1]$, one gets for some $C, \beta$,

(4.31) $\quad \|M_{h,1}^n\|_{L^\infty \to L^\infty} \leq C \int_{\widetilde{\tau}_1(h)}^\infty e^{-nh^2 x} (1+x)^\beta \, dx,$

and we get for some $C'$,

(4.32) $\quad \|M_{h,1}^n\|_{L^\infty \to L^\infty} \leq C' e^{-nh^2 \widetilde{\tau}_1(h)} \qquad \forall n \geq C_0 h^{-2}.$

Since $M_h^n$ is bounded by 1 on $L^\infty$, we get from $M_h^n - \Pi_0 = M_{h,1}^n + M_{h,2}^n$ and (4.31) that there exist $C_1, m$ such that $\|M_{h,2}^n\|_{L^\infty \to L^\infty} \leq C_1 h^{-m}$ for all $n \geq 1$. Next we use (1.15) to write $M_h = m_h + \mathcal{K}_h$ with

(4.33)
$$\|m_h\|_{L^\infty \to L^\infty} \leq \gamma < 1,$$
$$\|\mathcal{K}_h\|_{L^2 \to L^\infty} \leq C_2 h^{-d/2}.$$

From this, we deduce that for any $p = 1, 2, \ldots$ one has $M_h^p = A_{p,h} + B_{p,h}$, with $A_{1,h} = m_h, B_{1,h} = \mathcal{K}_h$ and the recurrence relation $A_{p+1,h} = m_h A_{p,h}, B_{p+1,h} = m_h B_{p,h} + \mathcal{K}_h M_h^p$. Thus one gets since $M_h^p$ is bounded by 1 on $L^2$,

(4.34)
$$\|A_{p,h}\|_{L^\infty \to L^\infty} \leq \gamma^p,$$
$$\|B_{p,h}\|_{L^2 \to L^\infty} \leq C_2 h^{-d/2}(1 + \gamma + \cdots + \gamma^p) \leq C_2 h^{-d/2}/(1-\gamma).$$

Observe that $\|M_{h,2}^n\|_{L^\infty \to L^2} \leq \|M_{h,2}^n\|_{L^2 \to L^2} \leq \delta^n$ and for $q, p \geq 1$, one gets, using (4.34),

(4.35)
$$\|M_{h,2}^{p+q}\|_{L^\infty \to L^\infty} = \|M_h^p M_{h,2}^q\|_{L^\infty \to L^\infty}$$
$$\leq \|A_{p,h} M_{h,2}^q\|_{L^\infty \to L^\infty} + \|B_{p,h} M_{h,2}^q\|_{L^\infty \to L^\infty}$$
$$\leq C_1 h^{-m} \gamma^p + C_2 h^{-d/2} \delta^q/(1-\gamma),$$

and this implies for some $C, \mu > 0$,

(4.36) $$\|M_{h,2}^n\|_{L^\infty \to L^\infty} \leq C e^{-n\mu} \qquad \forall n \geq 1/h,$$

and thus the contribution of $M_{h,2}^n$ is far smaller than the bound we have to prove in (4.28). The proof of Theorem 3 is complete.

## APPENDIX: CONVERGENCE TO THE BROWNIAN MOTION

The purpose of this appendix is to answer a question of one of the referees about the convergence of the previous Metropolis chain to the Brownian motion on a Riemannian manifold $(M, g)$. One classical and efficient way to prove such convergence is the use of Dirichlet forms (see [9]). Here, we present a self-contained proof, in the spirit of ([12], Chapter 2.4), making use of our previous results. The two main estimates are: the large deviation estimate (A.15) of Proposition 1, and the "central limit" theorem (A.46) of Proposition 2.

We refer to [8] and [11] for a construction of the Brownian motion on $(M, g)$. For a given $x_0 \in M$, let $X_{x_0} = \{\omega \in C^0([0, \infty[, M), \omega(0) = x_0\}$ be the set of continuous paths from $[0, \infty[$ to $M$, starting at $x_0$, equipped with the topology of uniform convergence on compact subsets of $[0, \infty[$, and let $\mathcal{B}$ be the Borel $\sigma$-field generated by the open sets in $X_{x_0}$. Let $W_{x_0}$ be the Wiener measure on $X_{x_0}$, and let $p_t(x, y) \, d_g y$ be the heat kernel, that is, the kernel of the self-adjoint operator $e^{t\Delta_g/2}$. Then $W_{x_0}$ is the unique probability on $(X_{x_0}, \mathcal{B})$, such that for any $0 < t_1 < t_2 < \cdots < t_k$ and any Borel sets $A_1, \ldots, A_k$ in $M$, one has

(A.1) $$W_{x_0}(\omega(t_1) \in A_1, \omega(t_2) \in A_2, \ldots, \omega(t_k) \in A_k)$$
$$= \int_{A_1 \times A_2 \times \cdots \times A_k} p_{t_k - t_{k-1}}(x_k, x_{k-1}) \cdots p_{t_2 - t_1}(x_2, x_1)$$
$$\times p_{t_1}(x_1, x_0) \, d_g x_1 \, d_g x_2 \cdots d_g x_k.$$

For $h \in \,]0, 1]$, let $\mathcal{M}_{h,x_0}^{\mathbb{N}}$ be the closed subset of the product space $M^{\mathbb{N}}$,

(A.2) $$\mathcal{M}_{h,x_0}^{\mathbb{N}} = \{\underline{x} = (x_1, x_2, \ldots, x_n, \ldots), \forall j \geq 0, d_g(x_j, x_{j+1}) \leq h\}.$$

Equipped with the product topology, $M^{\mathbb{N}}$ is a compact metrisable space, and the Metropolis chain starting at $x_0$ defines a probability $\mathcal{P}_{x_0,h}$ on $M^{\mathbb{N}}$, such that $\mathcal{P}_{x_0,h}(\mathcal{M}_{h,x_0}^{\mathbb{N}}) = 1$, by setting for all $k$ and all Borel sets $A_1, \ldots, A_k$ in $M$,

(A.3) $$\mathcal{P}_{x_0,h}(x_1 \in A_1, x_2 \in A_2, \ldots, x_k \in A_k)$$
$$= \int_{A_1 \times A_2 \times \cdots \times A_k} M_h(x_{k-1}, dx_k) \cdots M_h(x_1, dx_2) M_h(x_0, dx_1),$$

where the Metropolis kernel $M_h(x, dy)$ is defined in (1.15). Let $j_{x_0,h}$ be the map from $\mathcal{M}_{h,x_0}^{\mathbb{N}}$ into $X_{x_0}$ defined by

(A.4) $\qquad j_{x_0,h}(\underline{x}) = \omega \quad \Longleftrightarrow \quad \forall j \geq 0 \quad \omega(jh^2/(d+2)) = x_j$

and

(A.5) $\qquad \forall t \in \left[\dfrac{jh^2}{d+2}, \dfrac{(j+1)h^2}{d+2}\right] \quad \omega(t)$ is the geodesic curve connecting $x_j$ to $x_{j+1}$ at speed $h^{-2}(d+2)d_g(x_j, x_{j+1})$.

Observe that for $h > 0$ given, smaller than the injectivity radius of the Riemannian manifold $M$, the map $j_{x_0,h}$ is well defined and continuous. Let $P_{x_0,h}$ be the probability on $X_{x_0}$ defined as the image of $\mathcal{P}_{x_0,h}$ by the continuous map $j_{x_0,h}$. Our aim is to prove that $P_{x_0,h}$ converges weakly to the Wiener measure $W_{x_0}$ when $h \to 0$.

THEOREM 4. *For any bounded continuous function $\omega \mapsto f(\omega)$ on $X_{x_0}$, one has*

(A.6) $$\lim_{h \to 0} \int f \, dP_{x_0,h} = \int f \, dW_{x_0}.$$

Observe that the proof below shows that our study of the Metropolis chain on the manifold $M$ is also a way to prove the existence of the Brownian motion on $M$.

Let us recall that the Metropolis operator $M_h$ acting on $L^2 = L^2(M, d\mu_M)$ with $d\mu_M = d_g x / Vol(M))$ is defined by (1.17). If $\varphi$ is a Lipschitz function on $M$, we denote by $M_{h,\varphi}$ the bounded operator on $L^2$ defined by

$$M_{h,\varphi} = e^{\varphi/h} M_h e^{-\varphi/h}.$$

The first ingredient we use in the proof of Theorem 4 is the following lemma, which gives an $L^2$-estimate on the resolvent $(z - M_h)^{-1}$ near $z = 1$.

LEMMA 9. *Let $\psi$ be a real valued Lipschitz function on $M$, $\rho > 0$ and $0 < \theta < 2\pi$. Let us assume that the following inequality holds true:*

(A.7) $\qquad \rho \sin(\theta/2) - \displaystyle\sum_{k=2}^{\infty} \dfrac{\rho^{k/2} \|\psi\|_{\text{Lips}}^k}{k!} |\sin((k-1)\theta/2)| = c > 0.$

*Then, with $w = \rho e^{i\theta} \in \mathbb{C} \setminus [0, \infty[$ and $\varphi = i\rho^{1/2} e^{i\theta/2} \psi$, one has*

(A.8) $\qquad\qquad\qquad \|(1 - M_{h,\varphi} - w)^{-1}\|_{L^2} \leq 1/c.$

PROOF. If $k(x,y)$ is a complex valued bounded measurable function on $M \times M$, let $\mathcal{A}_{k,h}$ be the bounded operator on $L^2$,

$$(A.9) \quad \mathcal{A}_{k,h}(f)(x) = \int_{d_g(x,y) \leq h} \min\left(\frac{1}{|B(x,h)|}, \frac{1}{|B(y,h)|}\right) k(x,y) f(y) \, d_g y.$$

With $k^*(x,y) = \overline{k}(y,x)$, the adjoint on $L^2$ of $\mathcal{A}_{k,h}$ is equal to $\mathcal{A}_{k^*,h}$, and one has the obvious estimate

$$(A.10) \quad \|\mathcal{A}_{k,h}\|_{L^2} \leq \|k\|_{L^\infty(M \times M)}.$$

From (1.15) and (1.2), one has $M_h = m_h + \mathcal{A}_{1,h}$, and an easy calculation gives

$$(A.11) \quad M_{h,\varphi} = m_h + \mathcal{A}_{k_\varphi,h}, \qquad k_\varphi(x,y) = 1_{d_g(x,y) \leq h} e^{(\varphi(x) - \varphi(y))/h}.$$

Let $\tau(x,y) = 1_{d_g(x,y) \leq h} i(\psi(x) - \psi(y))/h$. With $\varphi = i\rho^{1/2} e^{i\theta/2} \psi$, we thus get

$$(A.12) \quad M_{h,\varphi} = m_h + \sum_{k=0}^\infty \frac{(\rho^{1/2} e^{i\theta/2})^k}{k!} \mathcal{A}_{\tau^k,h}.$$

From (A.12) and $w = \rho e^{i\theta}$, we get with $S = -e^{-i\theta/2}(1 - M_{h,\varphi} - w)$,

$$(A.13) \quad \begin{aligned} S &= -e^{-i\theta/2}(1 - M_h) + \rho^{1/2} \mathcal{A}_{\tau,h} + \rho e^{i\theta/2} Id + N, \\ N &= e^{-i\theta/2} \sum_{k=2}^\infty \frac{(\rho^{1/2} e^{i\theta/2})^k}{k!} \mathcal{A}_{\tau^k,h}. \end{aligned}$$

Since $\tau^* = \tau$, the second term in the first line of (A.13) is self-adjoint, and we get

$$(A.14) \quad \begin{aligned} \text{Im}(S) &= \sin(\theta/2)(1 - M_h) + \rho \sin(\theta/2) Id + \text{Im}(N), \\ \text{Im}(N) &= \sum_{k=2}^\infty \frac{\rho^{k/2}}{k!} \sin((k-1)\theta/2) \mathcal{A}_{\tau^k,h}. \end{aligned}$$

From $\sin(\theta/2)(1 - M_h) \geq 0$, and since from (A.10) the self-adjoint operator $\mathcal{A}_{\tau^k,h}$ has norm $\leq \|\psi\|_{\text{Lips}}^k$, we get from (A.7) and (A.14) that $\text{Im}(S) \geq cId$. The proof of Lemma 9 is complete. $\square$

From Lemma 9, we shall now deduce a key estimate on the probability that $X_{h,x_0}^n$, the $n$th step of the Metropolis chain starting at $x_0$, satisfies $d_g(X_{h,x_0}^n, x_0) > \varepsilon$. Let $\varepsilon_0 > 0$ be smaller than the injectivity radius of the Riemannian manifold $M$.

PROPOSITION 1. *There exist positive constants $C, A, a, c_0, h_0 > 0$ such that for all $\varepsilon \in ]0, \varepsilon_0]$, all $\delta \in ]0, c_0\varepsilon^2]$ and all $h \in ]0, h_0]$, the following inequality holds true:*

$$\text{(A.15)} \qquad \sup_{x_0 \in M, nh^2 \leq \delta} \mathcal{P}_{x_0,h}(d_g(X^n_{h,x_0}, x_0) > \varepsilon) \leq C\varepsilon^{-A} e^{-a\varepsilon^2/\delta}.$$

PROOF. We may assume $nh \geq \varepsilon$, since otherwise $\mathcal{P}_{x_0,h}(d_g(X^n_{h,x_0}, x_0) > \varepsilon) = 0$. In the proof, we denote by $a, A, C$ positive constants, changing from line to line, but which are independent of $h, \varepsilon$, $x_0 \in M$ and $n \geq 1$. One has

$$\text{(A.16)} \qquad \begin{aligned} \mathcal{P}_{x_0,h}(d_g(X^n_{h,x_0}, x_0) > \varepsilon) &= \int_{d_g(y,x_0)>\varepsilon} M_h^n(x_0, dy) \\ &= M_h^n(1_{d_g(y,x_0)>\varepsilon})(x_0). \end{aligned}$$

Let $\varphi(r) \in C^\infty([0, \infty[)$ be a nondecreasing function equal to 0 for $r \leq 3/4$ and equal to 1 for $r \geq 1$. For $\varepsilon \in ]0, \varepsilon_0]$ and $x_0 \in M$, set

$$\text{(A.17)} \qquad \varphi_{x_0,\varepsilon}(x) = \varphi\left(\frac{d_g(x, x_0)}{\varepsilon}\right).$$

Then $\varphi_{x_0,\varepsilon}$ is a smooth function, and from $1_{d_g(y,x_0)>\varepsilon} \leq \varphi_{x_0,\varepsilon} \leq 1$, we get, since $M_h$ is Markovian,

$$\text{(A.18)} \qquad M_h^n(1_{d_g(y,x_0)>\varepsilon}) \leq M_h^n(\varphi_{x_0,\varepsilon}) \leq M_h^n(1) = 1.$$

We first deduce from Lemma 9 the following estimates on $M_h^n(\varphi_{x_0,\varepsilon})$.

LEMMA 10. *There exists $c_0 > 0$ such that for $nh^2 \leq c_0\varepsilon^2$, the following inequalities hold true:*

$$\text{(A.19)} \qquad \|M_h^n(\varphi_{x_0,\varepsilon})\|_{L^2(B(x_0,\varepsilon/2))} \leq C e^{-a\varepsilon^2/nh^2};$$

$$\text{(A.20)} \qquad \|M_h^n(\varphi_{x_0,\varepsilon})\|_{L^\infty(B(x_0,\varepsilon/4))} \leq C h^{-d/2} e^{-a\varepsilon^2/nh^2}.$$

PROOF. By the Cauchy–Schwarz formula, the self-adjoint operator $M_h^n$ is equal to

$$\text{(A.21)} \qquad M_h^n = \frac{1}{2i\pi} \int_\sigma z^n (z - M_h)^{-1} dz,$$

where $\sigma$ is a contour in the complex plane surrounding the spectrum of $M_h$ with the counter-clockwise orientation. Let $\theta_0 \in ]0, \pi/2[$ close to $\pi/2$ and $\rho_0 > 0$ small be given. Since we know that the spectrum of $M_h$ is a subset of $[-\gamma, 1[$ with $\gamma \in [0, 1[$, we may choose $\sigma$ in the form $\sigma_1 \cup \sigma_2$, with

$$\sigma_1 = \{z = 1 - w(\theta), w(\theta) = \rho(\theta)e^{i\theta}, \theta \in [\theta_0, 2\pi - \theta_0]\},$$

where the function $\rho(\theta) > 0$ takes small values, is such that $\rho(\theta) = \rho(2\pi - \theta)$, $\rho_0 = \rho(\theta_0)$ and will be chosen later, and with $q = |1 - \rho_0 e^{i\theta_0}| < 1$,

$$\sigma_2 \subset \{|z| \leq q, \text{dist}(z, [-\gamma, 1]) \geq \rho_0 \sin(\theta_0)\}.$$

Set $g = \varphi_{x_0,\varepsilon}$ and $f_z = (z - M_h)^{-1}g$. For $z \in \sigma_2$, one has $\|f_z\|_{L^2} \leq \frac{\|g\|_{L^2}}{\rho_0 \sin(\theta_0)}$, and from $(z - m_h)f_z = \mathcal{A}_{1,h}f_z + g$, $|z - m_h(x)| \geq \text{dist}(z, [0,1]) \geq \rho_0 \sin(\theta_0)$, and $\|\mathcal{A}_{1,h}f_z\|_{L^\infty} \leq Ch^{-d/2}\|f_z\|_{L^2}$, we get for $z \in \sigma_2$, with a constant $C$ changing from line to line,

(A.22)
$$\|f_z\|_{L^\infty} \leq \frac{1}{\rho_0 \sin(\theta_0)}(\|\mathcal{A}_{1,h}f_z\|_{L^\infty} + \|g\|_{L^\infty})$$
$$\leq Ch^{-d/2}(\rho_0 \sin(\theta_0))^{-2}.$$

This gives

(A.23)
$$\left\| \int_{\sigma_2} z^n (z - M_h)^{-1}(g) \, dz \right\|_{L^2} \leq Cq^n (\rho_0 \sin(\theta_0))^{-1},$$
$$\left\| \int_{\sigma_2} z^n (z - M_h)^{-1}(g) \, dz \right\|_{L^\infty} \leq Cq^n h^{-d/2} (\rho_0 \sin(\theta_0))^{-2}.$$

Observe that since $nh \geq \varepsilon$, one has $q^n \leq e^{-a\varepsilon/h} \leq e^{-a\varepsilon^2/nh^2}$. Thus (A.23) gives (A.19) and (A.20) for the contribution of $\sigma_2$. Next we use Lemma 9 to bound the contribution of $\sigma_1$ in (A.21).

Let $\mu < 1$ and set $\psi(x) = \mu\sqrt{2}\,\text{dist}(x, B(x_0, \varepsilon/2))$. One has $\|\psi\|_{\text{Lips}} = \mu\sqrt{2}$, and if $\rho(\theta) > 0$ is small enough, inequality (A.7) is fulfilled with a constant $c \simeq \rho(\theta)\sin(\theta/2)(1 - \mu) + O(\rho^{3/2}(\theta)) \simeq \rho(\theta)$. From (A.8), and $(z - M_{h,\varphi})e^{\varphi/h}f_z = e^{\varphi/h}g$, we get for $z = 1 - w(\theta) \in \sigma_1$, since $\varphi = 0$ on $B(x_0, \varepsilon/2)$, $g = 0$ on $B(x_0, 3\varepsilon/4)$, and $|e^{\varphi/h}| = |e^{iw^{1/2}(\theta)\mu\sqrt{2}\,\text{dist}(x, B(x_0, \varepsilon/2))/h}|$,

(A.24)
$$\|f_z\|_{L^2(B(x_0,\varepsilon/2))} \leq \|e^{\varphi/h}f_z\|_{L^2} \leq \frac{C}{\rho(\theta)}\|e^{\varphi/h}g\|_{L^2}$$
$$\leq \frac{C'}{\rho(\theta)}|e^{iw^{1/2}(\theta)\mu\sqrt{2}\varepsilon/4h}|.$$

One has $(z - m_h)f_z = \mathcal{A}_{1,h}(f_z) + g$ with $g = 0$ on $B(x_0, \varepsilon/2)$, and $h \leq c_0\varepsilon$ since $h\varepsilon \leq nh^2 \leq c_0\varepsilon^2$. For $c_0 < 1/4$, we thus get from (A.24),

(A.25)   $\|f_z\|_{L^\infty(B(x_0,\varepsilon/4))} \leq C\rho^{-1}(\theta)h^{-d/2}|e^{iw^{1/2}a\varepsilon/h}|.$

On $\sigma_1$, we set $z = 1 - w = 1 - u^2$, $u = \rho^{1/2}(\theta)e^{i\theta/2} = w^{1/2}$. Then one has

(A.26)   $\int_{\sigma_1} z^n (z - M_h)^{-1}(g) \, dz = \int_\gamma (1 - u^2)^n f_{1-u^2} 2u \, du,$

where $\gamma$ is a contour in the upper half plane $\operatorname{Im}(u) > 0$ connecting $u_- = -\rho_0^{1/2} e^{-i\theta_0/2}$ to $u_+ = \rho_0^{1/2} e^{i\theta_0/2}$. From (A.24), (A.25) and (A.26), we deduce

(A.27)
$$\left\| \int_{\sigma_1} z^n (z - M_h)^{-1}(g) \, dz \right\|_{L^2(B(x_0, \varepsilon/2))} \leq CJ,$$
$$\left\| \int_{\sigma_1} z^n (z - M_h)^{-1}(g) \, dz \right\|_{L^\infty(B(x_0, \varepsilon/4))} \leq Ch^{-d/2} J,$$

where $J$ is defined by (with $a > 0$ small)

(A.28) $$J = \int_\gamma |(1 - u^2)^n e^{iua\varepsilon/h}| \frac{|du|}{|u|},$$

and it remains to verify that $J$ satisfies

(A.29) $$J \leq C_2 e^{-a_2 \varepsilon^2 / nh^2}.$$

At this point, we use the classical steepest descent method in order to choose the contour $\gamma$ such that (A.29) holds true. One has $(1-u^2)^n e^{iua\varepsilon/h} = e^{n(\log(1-u^2) + iru)}$ with $r = a\varepsilon/nh \in \,]0, a]$. Thus, $r > 0$ is a small parameter. The phase $\Phi(u) = \log(1 - u^2) + iru$ has a single nondegenerate critical point $u_c$ near 0, which satisfies, $u_c = ir/2 + \mathcal{O}(r^3)$, and the critical value is equal to $\Phi(u_c) = -r^2/4 + \mathcal{O}(r^4)$. Moreover, one has $\Phi''(u_c) = -2 + \mathcal{O}(r^2)$. It is then easy to verify that one can select the contour $\gamma$ in $\operatorname{Im}(u) \geq r/4$ connecting $u_-$ to $u_+$, and such that on $\gamma$, one has both $\operatorname{Re}(\Phi(u)) \leq \operatorname{Re}(\Phi(u_c)) - C_0|u - u_c|^2$ and $|u| \geq C_0(r + |u - u_c|)$ for some $C_0 > 0$. We thus get

(A.30) $$J \leq Ce^{n\Phi(u_c)} \int_{-\infty}^\infty e^{-ns^2} \frac{ds}{r + |s|}.$$

Then we get (A.29) from (A.30); one has $n\Phi(u_c) \leq -a^2\varepsilon^2/8nh^2$, and since $r\sqrt{n} = a\varepsilon/h\sqrt{n} \geq ac_0^{-1/2}$, one has $\int_{-\infty}^\infty e^{-ns^2} \frac{ds}{r+|s|} \leq C'/r\sqrt{n} \leq C_2$. The proof of Lemma 10 is complete. $\square$

Next, to deduce from the $L^2$ estimate (A.19) the desired $L^\infty$ estimate (A.15), we use the following lemma.

LEMMA 11. *For given $a_0, A_0, C_0$, there exist $a_1, A_1, C_1, p > 0, q > 0$ such that for $\varepsilon \in \,]0,1], n \geq 1$ and $0 < h \leq \varepsilon$, the following holds true: for any function $f$ on $M$ which satisfies $\|f\|_{L^\infty} \leq 1$, $\||\Delta_h|f\|_{L^\infty} \leq C_0 \varepsilon^{-2}$ and $\|f\|_{L^2(B(x_0, \varepsilon/2))} \leq C_0 \varepsilon^{-A_0} e^{-a_0 \varepsilon^2/nh^2}$, one has*

(A.31) $$\|f\|_{L^\infty(B(x_0, \varepsilon/4))} \leq C_1 (\varepsilon^{-A_1} e^{-a_1 \varepsilon^2/nh^2} + h^p \varepsilon^{-q}).$$

PROOF. Let $r_0 > 0$ and $\chi_0 \in C_0^\infty([0, 2r_0[)$ equal to 1 on $[0, r_0]$. Set $f_L = \chi_0(-h^2 \Delta_g) f$ and $f_H = f - f_L$. From Lemma 8, there exists $E_{1,h}$ and $\mathcal{N}_{1,h}$ such that $-f_H = E_{1,h}(1 - T_h) f + \mathcal{N}_{1,h} f$, and thus from (4.14) and $h^2|\Delta_h| = 2(d+2)(1-T_h)$, we get

$$\text{(A.32)} \qquad \|f_H\|_{L^\infty} \leq C h^2 \varepsilon^{-2}.$$

Let $\Phi_0 \in C_0^\infty([0, 2r_0[)$ be equal to 1 near the support of $\chi_0$. One has $\chi_0(1 - \Phi_0) = 0$ [we use the notation $\chi_0 = \chi_0(-h^2 \Delta_g), \Phi_0 = \Phi_0(-h^2 \Delta_g)$]. By Lemma 3 and with $|\Delta_h^0|$ defined by (3.31), we get

$$\text{(A.33)} \qquad \begin{aligned} \chi_0 |\Delta_h| f &= \chi_0 |\Delta_h^0| \Phi_0 f - 2(d+2) \chi_0 A_h f \\ &\quad - 2(d+2) \chi_0 (T_h h^{-2}(1 - \Phi_0)) f. \end{aligned}$$

Since $A_h \in \mathcal{E}_{cl}^{-\infty}$ and [by (A.32)] $\|h^{-2}(1-\Phi_0)f\|_{L^\infty} \leq C\varepsilon^{-2}$, we get

$$\text{(A.34)} \qquad \||\Delta_h^0| f_L\|_{L^\infty} = \|\chi_0 |\Delta_h^0| \Phi_0 f\|_{L^\infty} \leq C\varepsilon^{-2}.$$

By (1.6), one has $|\Delta_h^0| = -(1 + h^2 \Delta_g \widetilde{B}) \Delta_g$ with $\widetilde{B} \in \widetilde{\mathcal{E}}_{cl}^0$. Therefore, one has

$$|\Delta_h^0| f_L = |\Delta_h^0| \Phi_0 f_L = -(1 + h^2 \Delta_g \widetilde{B} \Phi_0) \Delta_g f_L.$$

If $r_0$ is small, the operator $1 + h^2 \Delta_g \widetilde{B} \Phi_0$ is invertible on $L^\infty$, and thus we get from (A.34)

$$\text{(A.35)} \qquad \|\Delta_g f_L\|_{L^\infty} \leq C\varepsilon^{-2}.$$

Let $\psi(x) \in [0,1]$ be a smooth function with support in the ball $B(x_0, \varepsilon/3)$ with $\psi(x)$ equal to 1 in the ball $B(x_0, \varepsilon/4)$, and such that $\|\nabla^\alpha \psi\|_{L^\infty} \leq C_\alpha \varepsilon^{-|\alpha|}$. Set $F(x) = \psi(x) f_L(x) = \psi(x)(f(x) - f_H(x))$. Using (A.32), $\Delta_g F = \psi \Delta_g f_L + [\Delta_g, \psi] f_L$ and (A.35), we get

$$\text{(A.36)} \qquad \begin{aligned} \|F\|_{L^2} &\leq C(\varepsilon^{-A_0} e^{-a_0 \varepsilon^2/nh^2} + h^2 \varepsilon^{-2+d/2}), \\ \|\Delta_g F\|_{L^\infty} &\leq C\varepsilon^{-2}, \qquad \|F\|_{L^\infty} \leq C. \end{aligned}$$

We now conclude that (A.31) holds true using (A.32), (A.36) and the classical interpolation inequality, with $\theta > \frac{d}{4+d}$

$$\text{(A.37)} \qquad \|F\|_{L^\infty} \leq C \|(1-\Delta_g) F\|_{L^\infty}^\theta \|F\|_{L^2}^{1-\theta}.$$

The proof of Lemma 11 is complete. $\square$

By the last inequality in (A.18) and (A.19), the function $f = M_h^n(\varphi_{x_0, \varepsilon})$ satisfies $\|f\|_{L^\infty} \leq 1$ and $\|f\|_{L^2(B(x_0, \varepsilon/2))} \leq C e^{-a\varepsilon^2/nh^2}$. Let us show that it satisfies also $\||\Delta_h| f\|_{L^\infty} \leq C\varepsilon^{-2}$. Let us recall that the operator $|\widetilde{\Delta}_h|$ is defined in (4.7). By (4.1) and (4.5), one has $|\Delta_h| = |\widetilde{\Delta}_h| + 2(d+2)h^{-2} R_h$ and

$\|R_h\|_{L^\infty} \leq Ch^3$. One gets easily from (2.17) $\||\Delta_h|\varphi_{x_0,\varepsilon}\|_{L^\infty} \leq C\varepsilon^{-2}$. Thus, one has also $\||\widetilde{\Delta}_h|\varphi_{x_0,\varepsilon}\|_{L^\infty} \leq C(\varepsilon^{-2}+h) \leq C'\varepsilon^{-2}$. Since $|\widetilde{\Delta}_h|$ commutes with $M_h$, one has $M_h^n(|\widetilde{\Delta}_h|\varphi_{x_0,\varepsilon}) = |\widetilde{\Delta}_h|M_h^n(\varphi_{x_0,\varepsilon})$, and this implies since $M_h$ is Markovian, $\||\widetilde{\Delta}_h|M_h^n(\varphi_{x_0,\varepsilon})\|_{L^\infty} \leq C\varepsilon^{-2}$. Thus we get $\||\Delta_h M_h^n(\varphi_{x_0,\varepsilon})\|_{L^\infty} \leq C(\varepsilon^{-2}+h) \leq C'\varepsilon^{-2}$. From Lemma 11, (A.16), (A.18) and (A.20) we thus get, for some $a, A, p, q > 0$,

$$
\begin{aligned}
\mathcal{P}_{x_0,h}(d_g(X_{x_0}^n, x_0) > \varepsilon) &\leq C(\varepsilon^{-A} e^{-a\varepsilon^2/nh^2} + h^p \varepsilon^{-q}), \\
\mathcal{P}_{x_0,h}(d_g(X_{x_0}^n, x_0) > \varepsilon) &\leq Ch^{-A} e^{-a\varepsilon^2/nh^2}.
\end{aligned}
\tag{A.38}
$$

Let $\alpha$ be such that $0 < \alpha < a/A$. It remains to observe that (A.38) implies (A.15), using the second line in case $h \geq e^{-\alpha\varepsilon^2/nh^2}$ and the first one if $h \leq e^{-\alpha\varepsilon^2/nh^2}$. The proof of Proposition 1 is complete. $\square$

With the result of Proposition 1, the proof of Theorem 4 follows now the classical proof of weak convergence of a sequence of random walks in the Euclidean space $\mathbb{R}^d$ to the Brownian motion on $\mathbb{R}^d$, for which we refer to ([12], Chapter 2.4). Let $T > 0$ be given. One has, for $0 < \delta \leq c_0\varepsilon^2$ and $h \in ]0, h_0]$,

$$
\begin{aligned}
\mathcal{P}_{x_0,h}(\exists j < l \leq h^{-2}T, &(l-j)h^2 \leq \delta, d_g(X_{x_0}^j, X_{x_0}^l) > 4\varepsilon) \\
&\leq \frac{C}{\delta} \sup_{y_0 \in M} \mathcal{P}_{y_0,h}(\exists j < l \leq h^{-2}\delta, d_g(X_{y_0}^j, X_{y_0}^l) > 4\varepsilon) \\
&\leq \frac{C}{\delta} \sup_{y_0 \in M} \mathcal{P}_{y_0,h}(\exists j \leq h^{-2}\delta, d_g(X_{y_0}^j, y_0) > 2\varepsilon) \\
&\leq \frac{2C}{\delta} \sup_{z_0 \in M, nh^2 \leq \delta} \mathcal{P}_{z_0,h}(d_g(X_{z_0}^n, z_0) > \varepsilon) \\
&\stackrel{\text{(by (A.15))}}{\leq} C'\delta^{-(1+A/2)} e^{-a\varepsilon^2/\delta}.
\end{aligned}
\tag{A.39}
$$

In fact, for the first inequality in (A.39), we just use the fact that the interval $[0,T]$ is a union of $\simeq C/\delta$ intervals of length $\delta/2$. The second inequality is obvious since the event $\{\exists j < l \leq h^{-2}\delta, d_g(X_{y_0}^j, X_{y_0}^l) > 4\varepsilon\}$ is a subset of $\{\exists j \leq h^{-2}\delta, d_g(X_{y_0}^j, y_0) > 2\varepsilon\}$. For the third, we use the fact that the event $A = \{\exists j \leq h^{-2}\delta, d_g(X_{y_0}^j, y_0) > 2\varepsilon\}$ is contained in $B \bigcup_{j<k}(C_j \cap D_j)$ with $B = \{d_g(X_{y_0}^k, y_0) > \varepsilon\}$ ($k$ is the greatest integer $\leq \delta h^{-2}$), $C_j = \{d_g(X_{y_0}^j, X_{y_0}^k) > \varepsilon\}$, $D_j = \{d_g(X_{y_0}^j, y_0) > 2\varepsilon \text{ and } d_g(X_{y_0}^l, y_0) \leq 2\varepsilon \text{ for } l < j\}$, and the fact that $C_j$ and $D_j$ are independent.

Using the definition (A.4), (A.5) of the map $j_{x_0,h}$, we get easily from (A.39) the convergence for $T > 0$ and $\varepsilon > 0$,

$$(A.40) \quad \lim_{\delta \to 0} \left( \limsup_{h \to 0} P_{x_0,h} \left( \max_{|s-t| \leq \delta, 0 \leq s,t \leq T} d_g(\omega(s), \omega(t)) > \varepsilon \right) \right) = 0.$$

Therefore, the family of probability $P_{x_0,h}$ is tight, hence is compact by the Prohorov theorem. It remains to verify that any weak limit $P_{x_0}$ of a sequence $P_{x_0,h_k}$, $h_k \to 0$, is equal to the Wiener measure $W_{x_0}$. By Theorem 4.15 of [12] we have to show that for any $m$, any $0 < t_1 < \cdots < t_m$, and any continuous function $f(x_1, \ldots, x_m)$, one has

$$(A.41) \quad \begin{aligned} &\lim_{k \to \infty} \int f(\omega(t_1), \ldots, \omega(t_m)) \, dP_{x_0,h_k} \\ &= \int f(x_1, \ldots, x_m) p_{t_m - t_{m-1}}(x_m, x_{m-1}) \cdots p_{t_2 - t_1}(x_2, x_1) \\ &\quad \times p_{t_1}(x_1, x_0) \, d_g x_1 \, d_g x_2 \cdots d_g x_m. \end{aligned}$$

As in [12], we may assume $m = 2$. For a given $t \geq 0$, let $n(t,h) \in \mathbb{N}$ be the greatest integer such that $h^2 n(t,h) \leq (d+2)t$. By (A.4), (A.5), one has $\operatorname{dist}(\omega(t), X_{h,x_0}^{n(t,h)}) \leq h$ and therefore $P_{x_0,h}(\operatorname{dist}(\omega(t), X_{h,x_0}^{n(t,h)}) > \varepsilon) = 0$ for $h \leq \varepsilon$. Thus we are reduced to prove

$$(A.42) \quad \begin{aligned} &\lim_{h \to 0} \int f(X_{h,x_0}^{n(t_1,h)}, X_{h,x_0}^{n(t_2,h)}) \, d\mathcal{P}_{x_0,h} \\ &= \int f(x_1, x_2) p_{t_2 - t_1}(x_2, x_1) p_{t_1}(x_1, x_0) \, d_g x_1 \, d_g x_2. \end{aligned}$$

From (A.3), one has

$$(A.43) \quad \begin{aligned} &\int f(X_{h,x_0}^{n(t_1,h)}, X_{h,x_0}^{n(t_2,h)}) \, d\mathcal{P}_{x_0,h} \\ &= \int f(x_1, x_2) M_h^{n(t_2,h) - n(t_1,h)}(x_1, dx_2) M_h^{n(t_1,h)}(x_0, dx_1). \end{aligned}$$

By (A.42), (A.43), we have to show that for any continuous function $f(x_1, x_2)$ on the product space $M \times M$, one has

$$(A.44) \quad \begin{aligned} &\lim_{h \to 0} \int_{M \times M} f(x_1, x_2) M_h^{n(t_2,h) - n(t_1,h)}(x_1, dx_2) M_h^{n(t_1,h)}(x_0, dx_1) \\ &= \int_{M \times M} f(x_1, x_2) p_{t_2 - t_1}(x_2, x_1) p_{t_1}(x_1, x_0) \, d_g x_1 \, d_g x_2, \end{aligned}$$

or, equivalently,

$$(A.45) \quad \begin{aligned} &\lim_{h \to 0} M_h^{n(t_1,h)} (M_h^{n(t_2,h) - n(t_1,h)} (f(x_1, \cdot))(x_1))(x_0) \\ &= e^{t_1 \Delta_g / 2} (e^{(t_2 - t_1) \Delta_g / 2} (f(x_1, \cdot))(x_1))(x_0). \end{aligned}$$

Since $\|M_h^{n(t,h)}\|_{L^\infty} \leq 1$ and $\|e^{t\Delta_g/2}\|_{L^\infty} \leq 1$, the following "central limit" theorem will conclude the proof of Theorem 4.

PROPOSITION 2. *For all $f \in C^0(M)$, and all $t > 0$, one has*

(A.46) $$\lim_{h \to 0} \|e^{t\Delta_g/2}(f) - M_h^{n(t,h)}(f)\|_{L^\infty} = 0.$$

PROOF. Since one has $\|M_h^{n(t,h)}\|_{L^\infty} \leq 1$ and $\|e^{t\Delta_g/2}\|_{L^\infty} \leq 1$, it is sufficient to prove that (A.46) holds true for $f \in \mathcal{D}$, with $\mathcal{D}$ a dense subset of the space $C^0(M)$, and therefore we may assume that $f = e_j$ is an eigenvector of $\Delta_g$. We set $n = n(t,h)$, and we use the notation of Section 4.2. From (4.36) and $n(t,h) \gg 1/h$, we get for some $a > 0$,

(A.47) $$\|M_{h,2}^{n(t,h)}(e_j)\|_{L^\infty} \leq Ce^{-at/h^2}.$$

One has

(A.48) $$\begin{aligned} g_h &= (M_{h,1}^{n(t,h)} + \Pi_0)e_j \\ &= \sum_{\widetilde{\tau}_k(h) \leq (1-\delta)h^{-2}} (1 - h^2\widetilde{\tau}_k(h))^{n(t,h)} \widetilde{e}_k^h(x) \int_M \widetilde{e}_k^h(y) e_j(y)\, d_g y. \end{aligned}$$

Let $A_j = \{k; |\widetilde{\tau}_k(h) - \frac{\lambda_j}{2(d+2)}| \leq \varepsilon\}$ with $\varepsilon$ small. Then from (4.8) and Theorem 2, one has $\sharp A_j = m_j = \dim Ker(\Delta_g + \lambda_j)$, and for any $k \notin A_j$, $|\int_M \widetilde{e}_k^h(y) \times e_j(y)\, d_g y| \leq C_k h$. Using (4.9), one has $|\widetilde{\tau}_k(h) - \frac{\lambda_k}{2(d+2)}| \leq C_k h$ for any given $k$. Take $N$ large and split the sum in (A.48) in the two pieces $\widetilde{\tau}_k(h) \leq N$ and $\widetilde{\tau}_k(h) > N$. Using the $L^\infty$ estimate (4.12) and the Weyl estimate (4.11) to bound the contribution of the sum on $\widetilde{\tau}_k(h) > N$, we get that there exists $C, a > 0$ and for all $N$, a constant $C(N)$ such that

(A.49) $$\|g_h - e^{-t\lambda_j/2}\Pi_{j,h}(e_j)\|_{L^\infty} \leq hC(N) + Ce^{-atN},$$

where $\Pi_{j,h}$ is the orthogonal projector on the vector space spanned by the $\widetilde{e}_k^h$ for $k \in A_j$. Let $\Pi_j$ be the orthogonal projector on $Ker(\Delta_g + \lambda_j)$. From (4.8) and Theorem 2, one has $\|\Pi_{j,h} - \Pi_j\|_{L^2} \leq C_j h$. From (4.24), one has $\|\widetilde{e}_k^h\|_{W^{1,p_*}} \leq C(1 + \widetilde{\tau}_k(h))^\alpha$ for some $p_* > d, \alpha > 0$. This implies $\|\Pi_{j,h} - \Pi_j\|_{L^2 \to W^{1,p_*}} \leq C_j$, and by interpolation $\|\Pi_{j,h} - \Pi_j\|_{L^2 \to L^\infty} \leq C_j h^\mu$ for some $\mu > 0$. Then (A.49) implies

(A.50) $$\|g_h - e^{-t\lambda_j/2}e_j\|_{L^\infty} \leq C_j h^\mu + hC(N) + Ce^{-atN}.$$

Clearly, (A.47) and (A.50) imply (A.46). The proof of Proposition 2 is complete. $\square$

**Acknowledgments.** We thank Persi Diaconis for numerous discussions and Erwann Aubry for his help on the calculus of differential geometry occurring in Section 2.

DÉPARTEMENT DE MATHÉMATIQUES
UNIVERSITÉ DE NICE SOPHIA-ANTIPOLIS
PARC VALROSE 06108, NICE CEDEX 02
FRANCE
E-MAILS: Gilles.Lebeau@unice.fr
Laurent.Michel@unice.fr